\documentclass[12pt]{amsart}
\usepackage{amsmath,amssymb,enumerate,cite,graphicx}

\usepackage[sc]{mathpazo}

\def\figdir{.}

\newcommand\BS\boldsymbol
\newcommand\dif{\,\mathrm{d}}

\newcommand\deriv[2]{\frac{\mathrm{d}#1}{\mathrm{d}#2}}
\newcommand\parderiv[2]{\frac{\partial #1}{\partial #2}}

\newcommand\picturehere[1]{\includegraphics[width=0.5\textwidth]{#1}}
\newcommand\largepicturehere[1]{\includegraphics[width=\textwidth]{#1}}

\numberwithin{equation}{section}

\begin{document}

\title[Suppression of chaos at slow variables by rapidly mixing fast
  dynamics]%
{Suppression of chaos at slow variables by rapidly mixing fast
dynamics through linear energy-preserving coupling}

\author{Rafail V. Abramov}
%\thanks {Department of Mathematics, Statistics and Computer Science,
%University of Illinois at Chicago, 851 S. Morgan st., Chicago, IL 60607
%(abramov@math.uic.edu).}
          %{Put the URL for your home page here if you have one}

          %Use \thanks statements for acknowledgements of grants and
          %support. They will appear below all the authors' addresses, so be
          %specific about which author is thanking whom:
          %\thanks{}

\address{Department of Mathematics, Statistics and Computer
  Science\\University of Illinois at Chicago\\851 S. Morgan st. (M/C
  249)\\ Chicago, IL 60607}

\email{abramov@math.uic.edu}

\subjclass[2000]{37M,37N,60G}

\date{\today}

\pagestyle{myheadings}
%\markboth{Suppression of chaos at slow variables}{Rafail V. Abramov}

\begin{abstract}
Chaotic multiscale dynamical systems are common in many areas of
science, one of the examples being the interaction of the
low-frequency dynamics in the atmosphere with the fast turbulent
weather dynamics. One of the key questions about chaotic multiscale
systems is how the fast dynamics affects chaos at the slow variables,
and, therefore, impacts uncertainty and predictability of the slow
dynamics. Here we demonstrate that the linear slow-fast coupling with
the total energy conservation property promotes the suppression of
chaos at the slow variables through the rapid mixing at the fast
variables, both theoretically and through numerical simulations. A
suitable mathematical framework is developed, connecting the slow
dynamics on the tangent subspaces to the infinite-time linear response
of the mean state to a constant external forcing at the fast
variables. Additionally, it is shown that the uncoupled dynamics for
the slow variables may remain chaotic while the complete multiscale
system loses chaos and becomes completely predictable at the slow
variables through increasing chaos and turbulence at the fast
variables. This result contradicts the common sense intuition, where,
naturally, one would think that coupling a slow weakly chaotic system
with another much faster and much stronger mixing system would result
in general increase of chaos at the slow variables.
\end{abstract}

\maketitle

\section{Introduction}
\label{sec:intro}

Dynamical systems, where the evolution of variables is separated
between two or more different time scales, are common in the
atmospheric/ocean science, physics, chemistry, molecular dynamics, and
many other areas. The structure of these systems is typically
characterized by the existence of a special subset of slow variables,
which evolve on a much longer time scale than the rest of the
variables. In particular, one can think of the low-frequency
variability models in the atmospheric science, where the slow
variables, usually the large scale empirical orthogonal functions
describing the large-scale slowly-varying patterns in the atmosphere
(such as the Arctic or North Atlantic oscillations, for example), are
coupled with small-scale fast processes, which are often very chaotic,
turbulent and unpredictable (with respect to the slow time scale, that
is). The inclusion of the fast processes often made the direct
numerical computation of the dynamics on the slow time scale
prohibitively expensive in real-world applications, which led to the
development of multiscale computational methods
\cite{ELiuVan,FatVan}. These methods make use of the averaging
formalism \cite{Pap,Van,Vol} to allow for the large time
discretization steps for the computation of the slow part of the
dynamics. Additionally, often the dynamics on the fast time scales are
not resolved, which requires suitable, often stochastic,
parameterizations for the fast dynamics which are derived from the
observed statistics of the fast variables (see, for example,
\cite{MajTimVan,MajTimVan2,MajTimVan3} and references therein).

One of the key questions about the behavior of multiscale dynamics,
which, however, does not seem to be extensively addressed in the
literature, is the effect of the rapidly mixing turbulent fast
dynamics on the chaotic properties of the slow variables. Formally,
the chaotic behavior of the whole multiscale system is determined by
its largest Lyapunov characteristic exponent, which often develops at
fast, strongly chaotic and rapidly mixing variables. However, the
largest Lyapunov exponent of the whole system does not in practice
characterize chaos at the slow variables, as it is often observed that
the behavior of the slow variables can be projected far beyond the
Lyapunov time scale of the fast variables. In fact, the averaging
formalism in \cite{Pap,Van,Vol} is based on the convergence of the
slow dynamics to a limit system which does not have any fast
variables, and the convergence occurs for finitely large time
intervals despite the fact that the the fast dynamics become
``infinitely fast'' in this limit, and the corresponding Lyapunov time
scale becomes infinitely short.

Here we study the chaotic behavior of slow variables by applying the
averaging formalism to the dynamics on the tangent subspaces of the
slow variables in a two-scale dynamical system. We demonstrate, both
analytically and through numerical simulations, that, if the linear
coupling between the slow and fast variables preserves the total
energy of the system, the averaged dynamics systematically suppresses
chaos at the slow variables. It is also speculated that the decreased
ratio of the autocorrelation time scale to the advection time scale at
the fast, turbulent and rapidly mixing variables could cause this
effect. Remarkably, it is observed that the full coupled system can
become completely predictable and non-chaotic at the slow variables,
while the uncoupled dynamics for the slow variables alone remain
chaotic. The theory is also systematically extended onto the
explicitly time-dependent and stochastically forced two-scale
dynamics.

The manuscript is organized as follows. In Section \ref{sec:lorenz} we
start with a numerical experiment with an appropriately rescaled
two-scale Lorenz model \cite{Abr6,CroVan,FatVan,Lor}, which
demonstrates that, for only slightly affected one-point statistics
such as the mean state and variance, the chaos at the slow variables
is suppressed by the increasing turbulent mixing at the fast variables
through the forcing parameter. In Sections
\ref{sec:twoscale}--\ref{sec:suppression} we develop a mathematical
framework for the slow averaged dynamics on the tangent subspaces.
For the linear energy-preserving deterministic coupling, it is shown
that chaos at the slow variables is suppressed wherever the
infinite-time linear response of the mean state of the fast dynamics
to a constant external forcing is positive definite. It is shown that
the positive-definiteness of the linear response is governed by the
statistical linear stability of the fast dynamics under sufficiently
small external perturbations, which, in turn, seems to be implemented
by the low ratio of the autocorrelation time to the advection time
scale at the fast variables. In Section \ref{sec:lorenz_revisited} we
go back to the Lorenz model and confirm the theoretical predictions by
computing the approximate infinite-time linear response through the
quasi-Gaussian fluctuation-dissipation theorem and observing how its
increasing positive-definiteness correlates with suppression of chaos
at the slow variables. In Section \ref{sec:suppression_uncoupled} we
show that it is possible for the uncoupled slow dynamics to remain
chaotic while the complete system with slow-fast dynamics and linear
energy-preserving coupling loses chaos and becomes predictable at the
slow scales by accelerating the turbulent mixing at the fast
variables. In Section \ref{sec:extension} we generalize the developed
framework onto the time-dependent and noisy dynamics, where the linear
energy-preserving coupling may also depend on the slow time. Section
\ref{sec:slow_fast_dep} gives a brief exposition of what happens when
the energy-preserving coupling is nonlinear in both the slow and fast
variables, and has the explicit fast time dependence. In Section
\ref{sec:invariance} we show that the developed framework is also
applicable to the multiscale dynamics without the explicit scaling
parameter $\varepsilon$. Section \ref{sec:summary} summarizes the
results of this work, as well as outlines the future directions of
research.

\section{A remarkable behavior of the two-scale Lorenz model}
\label{sec:lorenz}

Here we consider a two-scale system of autonomous ordinary
differential equations of the form
\begin{equation}
\label{eq:dyn_sys}
\deriv{\BS x}t=\BS F(\BS x,\BS y),\qquad
\deriv{\BS y}t=\frac 1\varepsilon\BS G(\BS x,\BS y),
\end{equation}
where $\BS x=\BS x(t)\in\mathbb R^{N_x}$ are the slow variables, $\BS
y=\BS y(t)\in\mathbb R^{N_y}$ are the fast variables, $\BS F$ and $\BS
G$ are $N_x$ and $N_y$ vector-valued functions of $\BS x$ and $\BS y$,
respectively, and $\varepsilon\ll 1$ is a time-scale separation
parameter.

For the two-scale system in \eqref{eq:dyn_sys}, one can write the
averaged dynamics for $\BS x$ alone as $\varepsilon\to 0$ using the
formalism from \cite{Abr6,Pap,Van,Vol}. The averaging formalism
produces the averaged system for $\BS x$ in the form
\begin{equation}
\label{eq:dyn_sys_slow_limiting_x}
\deriv{\BS x}t=\langle\BS F\rangle(\BS x),\qquad
\langle\BS F\rangle(\BS x)=\int_{\mathcal A_{\BS x}}\BS F(\BS x,\BS z)
\dif\mu_{\BS x}(\BS z),
\end{equation}
where $\mu_{\BS x}$ is the invariant probability measure on the
attractor $\mathcal A_{\BS x}$ of the limiting fast dynamics given by
\begin{equation}
\label{eq:dyn_sys_fast_limiting_z}
\deriv{\BS z}s=\BS G(\BS x,\BS z),
\end{equation}
where $\BS z=\BS z(s)$, while $\BS x$ is a fixed constant parameter
for \eqref{eq:dyn_sys_fast_limiting_z} and, consequently, $\mu_{\BS
  x}$. Here we tacitly assume the ergodicity of $\mu_{\BS x}$, as well
as that \eqref{eq:dyn_sys_slow_limiting_x} constitutes an
approximation to \eqref{eq:dyn_sys} for the finite times in the limit
as $\varepsilon\to 0$ (for details, see \cite{Pap,Van,Vol} and
references therein). For the above formalism to work, the right-hand
side of the averaged system in \eqref{eq:dyn_sys_slow_limiting_x} must
not be $O(\varepsilon)$, otherwise, higher-order homogenization
techniques are needed \cite{Van}.

We choose the two-scale forced damped Lorenz model
\cite{Abr5,Abr6,CroVan,FatVan,Lor} for the computational study of the
dynamical properties of a two-scale slow-fast process with generic
features of atmospheric dynamics, such as the presence of linearly
unstable waves, strong nonlinearity, forcing, dissipation, chaos and
mixing. The two-scale forced damped Lorenz model is given by
\begin{subequations}
\label{eq:lorenz_two_scale}
\begin{equation}
\dot x_i=x_{i-1}(x_{i+1}-x_{i-2})-x_i+F_x-\frac{\lambda_y}J
\sum_{j=1}^Jy_{i,j},
\end{equation}
\begin{equation}
\dot y_{i,j}=\frac 1\varepsilon\left[y_{i,j+1}
(y_{i,j-1}-y_{i,j+2})-y_{i,j}+F_y+\lambda_xx_i\right],
\end{equation}
\end{subequations}
where $1\leq i\leq N_x$, $1\leq j\leq J$. The following notations are
adopted above:
\begin{itemize}
\item $\BS x$ is the set of the slow variables of size $N_x$. The
  following periodic boundary conditions hold for $\BS x$:
  $x_{i+N_x}=x_i$;
\item $\BS y$ is the set of the fast variables of size $N_y=N_xJ$
  where $J$ is a positive integer. The following boundary conditions
  hold for $\BS y$: $y_{i+N_x,j}=y_{i,j}$ and $y_{i,j+J}=y_{i+1,j}$;
\item $F_x$ and $F_y$ are the constant forcing parameters;
\item $\lambda_x$ and $\lambda_y$ are the coupling parameters;
\item $\varepsilon$ is the time scale separation parameter.
\end{itemize}
Originally in \cite{CroVan,FatVan,Lor} there was no constant forcing
$F_y$ term in the equation for $\BS y$-variables in
\eqref{eq:lorenz_two_scale}, however, in its absence the behavior of
the $\BS y$-variables is strongly dissipative \cite{Abr5,Abr6}. Here,
as in \cite{Abr6}, we add a constant forcing $F_y$ in the right-hand
side of the second equation in \eqref{eq:lorenz_two_scale} to induce
the strongly chaotic behavior of the $\BS y$-variables with large
positive Lyapunov exponents. It is demonstrated in \cite{FatVan} that
the averaging formalism for \eqref{eq:lorenz_two_scale} is valid.

In the Lorenz model \eqref{eq:lorenz_two_scale}, $F_x$ and $F_y$
regulate the chaos and mixing of the $\BS x$ and $\BS y$ variables,
respectively \cite{Abr5,Abr6,Abr7,AbrMaj,AbrMaj4,AbrMaj5,MajAbrGro}.
However, the mean state and mean energy are also affected by the
changes in forcing, which affects the mean and energy trends in
coupling for the fixed coupling parameters. To adjust the effect of
coupling independently of forcing, here we perform the rescaling of
the Lorenz model as in \cite{MajAbrGro}. Consider the uncoupled Lorenz
model
\begin{equation}
\label{eq:L96}
\deriv{}tx_i=x_{i-1}(x_{i+1}-x_{i-2})-x_i+F
\end{equation}
with the same periodic boundary conditions as above \cite{LorEma}.
Observe that the long term statistical mean state $\bar x$ and the
standard deviation $\beta$ in \eqref{eq:L96} are the same for all
$x_i$ due to the translational invariance. Now, we rescale $\BS x$ and
$t$ as
\begin{equation}
x_i=\bar x+\beta q_i,\quad t=\frac\tau\beta,
\end{equation}
where the new variables $\BS q$ have zero mean state and unit standard
deviation. In the rescaled variables, the Lorenz model becomes
\begin{equation}
\label{eq:L96_rescaled}
\deriv{}\tau q_i=q_{i-1}(q_{i+1}-q_{i-2})+\frac 1\beta\left[\bar x(q_{i+1}
-q_{i-2})-q_i\right]+\frac{F-\bar x}{\beta^2},
\end{equation}
where $\bar x$ and $\beta$ are, of course, the functions of $F$. In
addition to setting the mean state and variance of $q_i$ to zero and
one, respectively, due to the time rescaling the autocorrelation
functions of $\BS z$ acquire identical time scaling for any $F$ (for
details, see \cite{MajAbrGro}).

Here, we similarly rescale the two-scale Lorenz model from
\eqref{eq:lorenz_two_scale}:
\begin{subequations}
\label{eq:lorenz_rescaled}
\begin{equation}
\dot x_i=x_{i-1}(x_{i+1}-x_{i-2})+\frac 1{\beta_x}
\left(\bar x(x_{i+1}-x_{i-2})-x_i\right)+\frac{F_x-\bar x}
{\beta_x^2}-\frac{\lambda_y}J\sum_{j=1}^Jy_{i,j},
\end{equation}
\begin{equation}
\dot y_{i,j}=\frac 1\varepsilon\left[y_{i,j+1}
(y_{i,j-1}-y_{i,j+2})+\frac 1{\beta_y}
\left(\bar y(y_{i,j-1}-y_{i,j+2})-y_{i,j}\right)+\frac{F_y-\bar y}
{\beta_y^2}+\lambda_xx_i\right],
\end{equation}
\end{subequations}
where $\bar x$, $\bar y$, $\beta_x$ and $\beta_y$ are the long term
means and variances of the corresponding uncoupled system in
\eqref{eq:L96} with either $F_x$ or $F_y$ set as a constant
forcing. Additionally, we write the fast limiting dynamics for
\eqref{eq:lorenz_rescaled} as
\begin{equation}
\label{eq:lorenz_fast_limiting}
\dot z_{i,j}=z_{i,j+1}
(z_{i,j-1}-z_{i,j+2})+\frac 1{\beta_y}\left(\bar z(z_{i,j-1}-z_{i,j+2})
-z_{i,j}\right)+\frac{F_y-\bar z}{\beta_y^2}+\lambda_xx_i,
\end{equation}
where $\BS x$ is given as an external parameter, as in
\eqref{eq:dyn_sys_fast_limiting_z}. The rescaling of the model implies
that the average mean states of $\BS x$ and $\BS y$ are close to zero,
and in this case the slow averaged system simply becomes the rescaled
Lorenz model in \eqref{eq:L96_rescaled}, where the right-hand side is
$O(1)$. Therefore, the homogenization techniques from \cite{Van} are
not needed, just like for the unrescaled Lorenz model in
\cite{FatVan}.

Below, we show the computed statistics of the rescaled Lorenz model
\eqref{eq:lorenz_rescaled} with the following parameters: $N_x=10$,
$N_y=40$, $F_x=6$, $F_y=6,8,12,16$ and $24$,
$\lambda_x=\lambda_y=0.25$, $\varepsilon=0.01$ (the time scale
separation between $\BS x$ and $\BS y$ is 100 times). The slow forcing
parameter $F_x=6$ is chosen so that the slow dynamics are not too
chaotic, mimicking the behavior of low-frequency variability in the
atmosphere (it is known from the previous work
\cite{Abr5,Abr6,Abr7,AbrMaj,AbrMaj4,AbrMaj5,MajAbrGro} that for $F=6$
the dynamics of the uncoupled model in \eqref{eq:L96} are weakly
chaotic). The coupling parameters $\lambda_x$ and $\lambda_y$ are set
to $0.25$ so that they are neither too small, nor too large, to ensure
the rich interaction between the slow and fast variables without
linearizing the rescaled Lorenz system too much. The time-scale
separation parameter $\varepsilon=0.01$ is, again, chosen so that it
is neither too large, nor too small (the time scale separation by two
orders of magnitude is consistent, for instance, with the difference
between annual and diurnal cycles in the atmosphere). Here and
everywhere else, the simple 2nd order Runge-Kutta scheme is used to
integrate the model in \eqref{eq:lorenz_rescaled} forward in time with
the discrete time step $\Delta t=10^{-5}$ and the time averaging
window of 10000 time units.

In the rescaled Lorenz model \eqref{eq:lorenz_rescaled}, it turns out
that the values of $F_x$ and $F_y$ do not significantly affect the
mean state and mean energy for both the slow variables $\BS x$ and
fast variables $\BS y$. To show this, in Table \ref{tab:mean_var_1} we
display the mean states and variances of both $\BS x$ and $\BS y$ for
the rescaled Lorenz model in \eqref{eq:lorenz_rescaled}.
\begin{table}%
\begin{center}%
\begin{tabular}{|c|}%
\hline
$N_x=10$, $N_y=40$, $F_x=6$, $\lambda_x=\lambda_y=0.25$, $\varepsilon=0.01$ \\
\hline
\begin{tabular}{c||c|c|c|c}
$F_y$ & $\BS x$-mean & $\BS x$-var & $\BS y$-mean & $\BS y$-var \\
\hline\hline
$6$ & $9.64\cdot 10^{-3}$ & $0.9451$ & $-2.38\cdot 10^{-3}$ & $1.066$ \\
$8$ & $2.817\cdot 10^{-2}$ & $0.9514$ & $-1.466\cdot 10^{-2}$ & $1.098$ \\
$12$ & $2.05\cdot 10^{-2}$ & $0.9336$ & $-2.719\cdot 10^{-2}$ & $1.139$ \\
$16$ & $-1.353\cdot 10^{-2}$ & $0.9006$ & $-4.028\cdot 10^{-2}$ & $1.153$ \\
$24$ & $-6.972\cdot 10^{-2}$ & $0.8434$ & $-6.075\cdot 10^{-2}$ & $1.167$ \\
\end{tabular}\\
\hline
\end{tabular}%
\end{center}%
\caption{The mean states and variances of the $\BS x$ and $\BS y$
  variables for the rescaled Lorenz model in
  \eqref{eq:lorenz_rescaled} with the following parameters: $N_x=10$,
  $N_y=40$, $F_x=6$, $F_y=6,8,12,16$ and $24$,
  $\lambda_x=\lambda_y=0.25$, $\varepsilon=0.01$.}
\label{tab:mean_var_1}
\end{table}
Observe that, despite different forcing regimes, the means and
variances for both $\BS x$ and $\BS y$ are almost unchanged, the mean
states being near zero while the variances being near one, as designed
by the rescaling for the uncoupled model in
\eqref{eq:L96_rescaled}. Here note that while the rescaling was
carried out for the corresponding uncoupled model (where it sets the
mean state to zero and variance to one precisely), using the same
rescaling parameters in the coupled model \eqref{eq:lorenz_rescaled}
still sets its means and variances near prescribed values zero and
one, respectively (although not precisely).

\begin{figure}%
\largepicturehere{\figdir/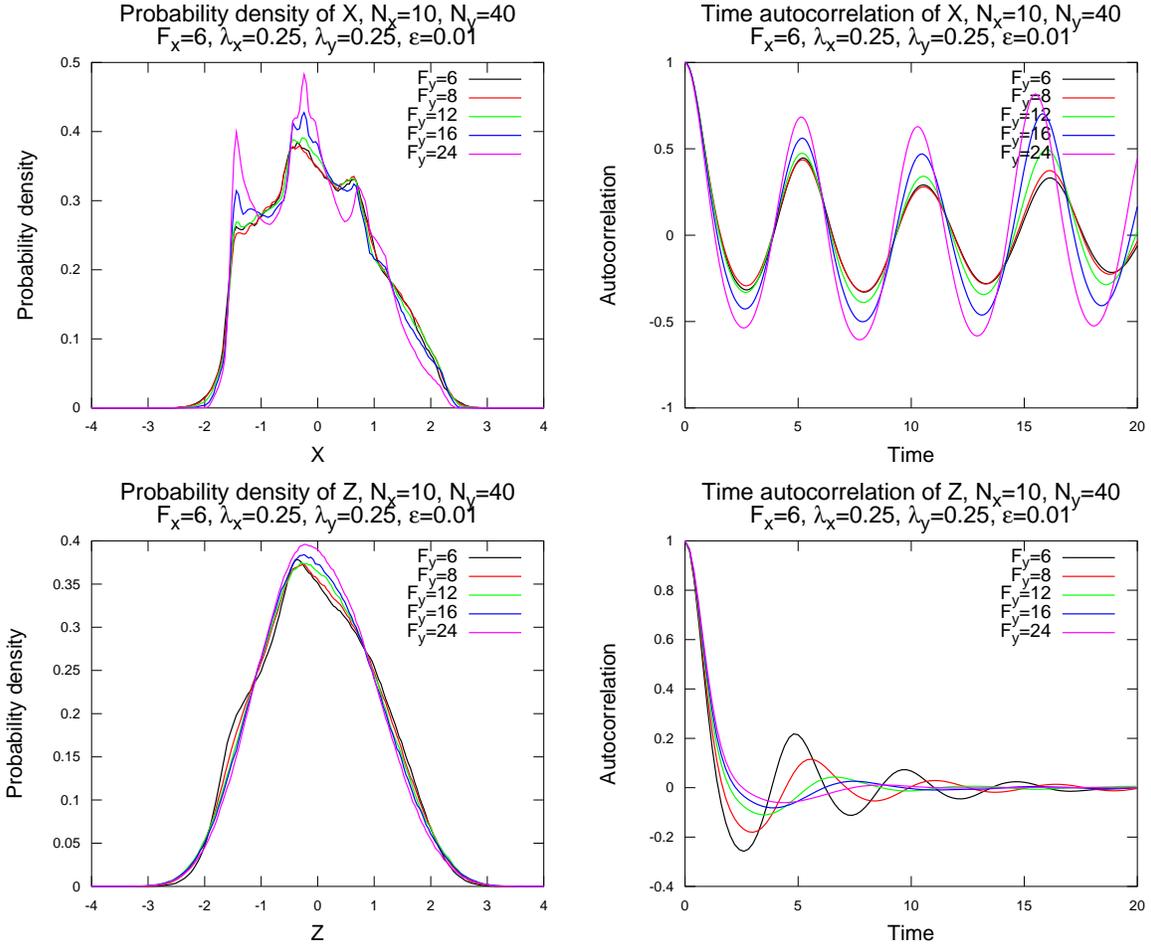}%
\caption{Probability density functions and time autocorrelation
  functions for $\BS x$ of \eqref{eq:lorenz_rescaled} and $\BS z$ of
  \eqref{eq:lorenz_fast_limiting} (with free parameter $\BS x$ set to
  its mean value) with the following parameters: $N_x=10$, $N_y=40$,
  $F_x=6$, $F_y=6,8,12,16$ and $24$, $\lambda_x=\lambda_y=0.25$,
  $\varepsilon=0.01$.}%
\label{fig:pdf_corr_1}%
\end{figure}

However, the different values of the fast forcing $F_y$ introduce some
changes in the statistics (such as the probability density functions
and the time autocorrelation functions) of both the slow variables
$\BS x$ of \eqref{eq:lorenz_rescaled}, and the fast limiting variables
$\BS z$ of \eqref{eq:lorenz_fast_limiting}. To show this, in Figure
\ref{fig:pdf_corr_1} we demonstrate the probability density functions
and the time autocorrelation functions of both $\BS x$ and $\BS z$
(where the free parameter $\BS x$ is set to its mean value).
Observe that in Figure \ref{fig:pdf_corr_1} different values of the
fast forcing $F_y$ do not significantly affect the probability density
functions for both $\BS x$ and $\BS z$ -- they are quite similar with
some variations. Perhaps, the most noticeable qualitative change can
be observed in the time autocorrelation functions for both $\BS x$ and
$\BS z$ -- for smaller values of $F_y$ the time autocorrelation
functions of $\BS x$ oscillate with somewhat smaller amplitude (which
indicates somewhat stronger mixing), while the time autocorrelation
functions of $\BS z$ display visibly longer decay, which indicates
weaker mixing.

\begin{figure}%
\begin{center}%
\picturehere{\figdir/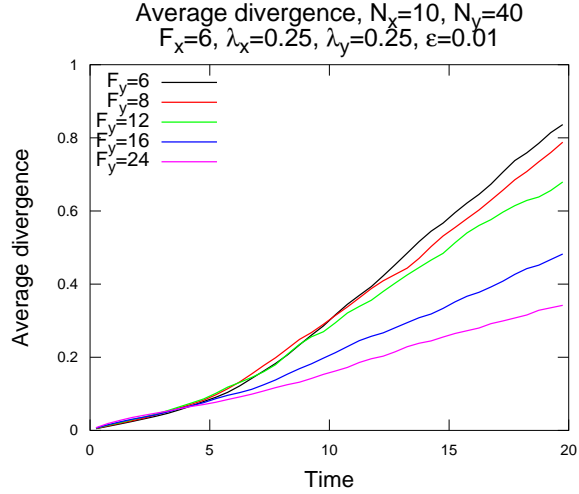}%
\end{center}%
\caption{Average divergence between perturbed and unperturbed running
  averages of the slow variables for the rescaled Lorenz model in
  \eqref{eq:lorenz_rescaled} with the following parameters: $N_x=10$,
  $N_y=40$, $F_x=6$, $F_y=6,8,12,16$ and $24$,
  $\lambda_x=\lambda_y=0.25$, $\varepsilon=0.01$.}%
\label{fig:av_diff_1}%
\end{figure}
\begin{figure}%
\picturehere{\figdir/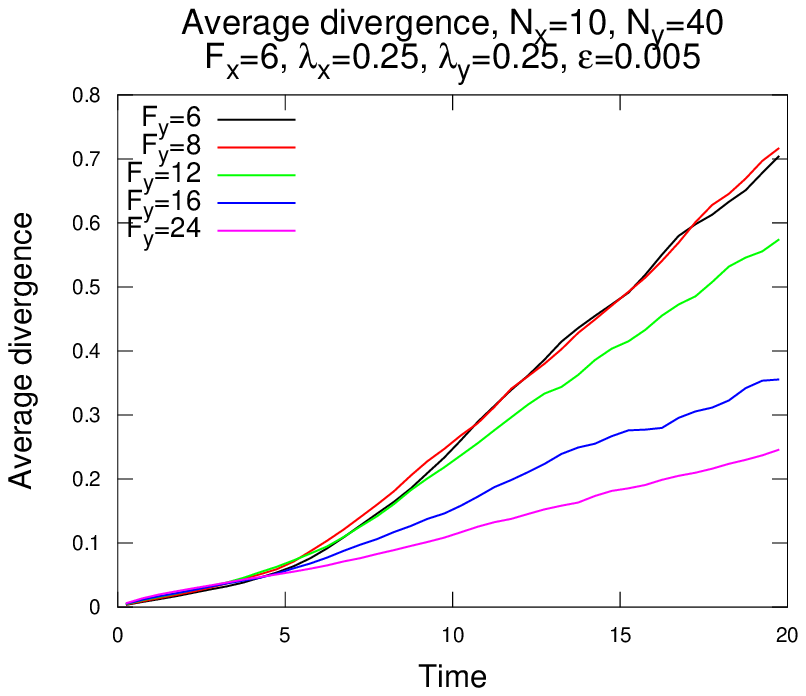}%
\picturehere{\figdir/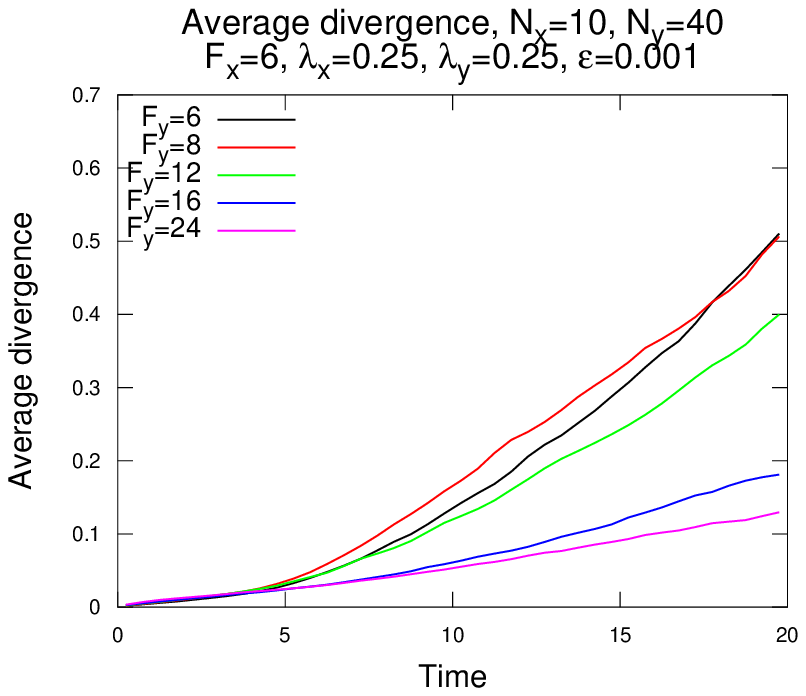}%
\caption{Average divergence between perturbed and unperturbed running
  averages of the slow variables for the rescaled Lorenz model in
  \eqref{eq:lorenz_rescaled} with the following parameters: $N_x=10$,
  $N_y=40$, $F_x=6$, $F_y=6,8,12,16$ and $24$,
  $\lambda_x=\lambda_y=0.25$, $\varepsilon=0.005$ and
  $\varepsilon=0.001$.}%
\label{fig:av_diff_2}%
\end{figure}

%
Having observed the statistics, we turn our attention to the chaotic
behavior of the slow variables $\BS x$ for the same range of
parameters. As the slow variables $\BS x$ are coupled to the fast
variables $\BS y$, the straightforward computation of the largest
Lyapunov exponent will characterize the chaos at the fast variables
$\BS y$, and, therefore, a different method should be used to quantify
chaos at the slow variables. Here we are going to observe the average
divergence behavior in time between the short-time (half of the time
unit) running averages $\langle\BS x\rangle(t)$ of the slow time
series $\BS x(t)$, which are initially generated very closely to each
other. Namely, we perform the following computation:
\begin{itemize}
\item We start with a generic initial condition $(\BS x_0,\BS y_0)$, and
  propagate it for 1000 time units to let it settle on the chaotic
  attractor of the rescaled Lorenz model in \eqref{eq:lorenz_rescaled};
\item A perturbed $(\BS x,\BS y)$-state is generated by a small random
  deviation (of order $\varepsilon$) from the computed
  trajectory. Then, the original trajectory $\BS x(t)$ and the
  perturbed state $\BS x^\prime(t)$ are integrated forward in time in
  parallel. The time evolution of the divergence $\|\langle\BS
  x^\prime\rangle(t)-\langle\BS x\rangle(t)\|$ between the original
  and the perturbed trajectory is recorded;
\item The latter operation is repeated 500 times with different
  snapshots of the same long-term trajectory, spaced by 20 time units,
  and averaged. The result is the time evolution of the divergence in
  $\langle\BS x\rangle(t)$ between the original an perturbed
  trajectory, averaged over the attractor (averaging window of 10000
  time units) of the rescaled Lorenz model in
  \eqref{eq:lorenz_rescaled}.
\end{itemize}
The short time-averaging window of half of the time unit for the
running average $\langle\BS x\rangle(t)$ ensures that the slow
variables $\BS x(t)$ do not change much during this window, while the
fast time series $\BS y(t)$ mix completely during the same short time
averaging window (see the time scales of the time autocorrelation
functions of $\BS x$ and $\BS y$ in Figure \ref{fig:pdf_corr_1} for
comparison). The results of this simulation for the same set of
parameters are shown in Figure \ref{fig:av_diff_1}. Remarkably, the
chaos in the slow $\BS x$-variables is consistently suppressed as the
fast forcing $F_y$ increases, as the unperturbed and perturbed slow
running averages $\langle\BS x\rangle(t)$ diverge from each other
slower and slower in time. It cannot be caused by the changing
statistical mean or variance of the slow or fast variables creating
average counteracting forcing at the slow variables, as Table
\ref{tab:mean_var_1} clearly indicates that the mean states and
variances of both the slow and fast variables do not change by a
significant amount for different $F_y$. This chaos suppression effect
at the slow variables appears to become somewhat stronger as the
time-scale separation is increased: in Figure \ref{fig:av_diff_2} we
show the same average divergence as in Figure \ref{fig:av_diff_1}, but
with the time-scale separation parameter $\varepsilon$ set to $0.005$
and $0.001$, respectively.

In what follows we set out to explain the observed effect of the chaos
suppression at the slow variables by increasing the turbulent mixing
at the fast variables, without affecting the mean forcing or
variability. As the effect pertains when $\varepsilon$ is decreased,
one would expect the averaging formalism \cite{Pap,Van,Vol} in the
limit as $\varepsilon\to 0$ to provide an approximate way to study the
effect from the dynamical system perspective, which we systematically
develop below.

\section{A two-scale slow-fast system and its tangent linear model for
the averaged slow variables}
\label{sec:twoscale}

To explain the observations for the two-scale rescaled Lorenz model in
the previous section, here we study the behavior of the tangent linear
model of the averaged dynamics in \eqref{eq:dyn_sys_slow_limiting_x}.
Below, we do not attempt at rigorous proofs of stated claims, instead
making reasonable assumptions about the system under consideration,
and justifying the claims under the assumptions made.

We start by deriving the tangent linear model for
\eqref{eq:dyn_sys_slow_limiting_x}. Here, the necessary assumption is
that the flow (or, as it is also called, the solution operator) $\BS
x(t)=\BS X(\BS x,t)$, generated by \eqref{eq:dyn_sys_slow_limiting_x},
is differentiable with respect to its initial condition $\BS x$. In
other words, we assume that the averaged tangent map
\begin{equation}
\BS{TX}^t_{\BS x}=\parderiv{}{\BS x}\BS X(\BS x,t)
\end{equation}
for \eqref{eq:dyn_sys_slow_limiting_x} exists and is given by the
equation
\begin{equation}
\label{eq:dyn_sys_slow_tangent_1}
\parderiv{}t\BS{TX}^t_{\BS x}=\parderiv{\langle\BS F\rangle}
{\BS x}(\BS x)\BS{TX}^t_{\BS x}.
\end{equation}
The partial derivative in \eqref{eq:dyn_sys_slow_tangent_1} is given,
formally, as
\begin{equation}
\label{eq:F_average}
\parderiv{\langle\BS F\rangle}{\BS x}(\BS x)=\int_{\mathcal A_{\BS
    x}}\parderiv{\BS F}{\BS x}(\BS x,\BS z) \dif\mu_{\BS x}(\BS
z)+\int_{\mathcal A_{\BS x}} \BS F(\BS x,\BS z)\dif\mu^\prime_{\BS
  x}(\BS z).
\end{equation}
Above we assume that the invariant measure $\mu_{\BS x}$
differentiably depends on $\BS x$, that is, for any $\mu_{\BS
  x}$-measurable $f(\BS z)$ we have
\begin{equation}
\label{eq:temp0}
\begin{split}
\int_{\mathcal A_{\BS x}} f(\BS z)\dif\mu_{\BS x+\delta\BS x}
(\BS z)-\int_{\mathcal A_{\BS x}} f(\BS z)\dif\mu_{\BS x}(\BS z)=\\=
\left(\int_{\mathcal A_{\BS x}} f(\BS z)\dif\mu_{\BS x}^\prime(\BS z)\right)
\delta\BS x+o(\|\delta\BS x\|),
\end{split}
\end{equation}
where $\mu_{\BS x}^\prime$ denotes the derivative of $\mu_{\BS x}$
with respect to $\BS x$. The differentiability property above requires
structural stability of $\mu_{\BS x}$ under changes in $\BS x$. It is
known that uniformly hyperbolic diffeomorphisms on the whole $\mathbb
R^{N_y}$ (Anosov) or just the nonwandering set of $\mathbb R^{N_y}$
(Axiom A) are structurally stable \cite{EckRue,Rue,Rue2,Rue3,You}.

Using the linear response theory
\cite{Abr5,Abr6,Abr7,AbrMaj4,AbrMaj5,AbrMaj6,MajAbrGro,Rue2,Ris}, the
second term in \eqref{eq:temp0} can be computed as the infinite-time
linear response operator of \eqref{eq:dyn_sys_fast_limiting_z} for a
perturbation of the external parameter $\BS x$. The calculations are
given in the Appendix \ref{sec:inf_time_lin_resp}, and the result is
\begin{equation}
\begin{split}
\int_{\mathcal A_{\BS x}}\BS F(\BS x,\BS z)\dif\mu_{\BS x}^\prime(\BS
z)=\int_0^\infty\int_{\mathcal A_{\BS x}}\parderiv{\BS F}{\BS y} (\BS
x,\phi_{\BS x}^s\BS z) \BS T_{\BS x,\BS z}^s\parderiv{\BS G}{\BS
  x}(\BS x, \BS z)\dif\mu_{\BS x}(\BS z)\dif s,
\end{split}
\end{equation}
where $\phi_{\BS x}^t$ is the flow, generated by
\eqref{eq:dyn_sys_fast_limiting_z}. For convenience, further we denote
the left-hand side of the above formula as $\BS H(\BS x)$. Thus, for
the averaged tangent dynamics in \eqref{eq:dyn_sys_slow_tangent_1} we
obtain
\begin{subequations}
\label{eq:dyn_sys_slow_tangent}
\begin{equation}
\parderiv{}t\BS{TX}^t_{\BS x}=\left(\left\langle\parderiv{\BS F}{\BS x}
\right\rangle(\BS x)+\BS H(\BS x)\right)\BS{TX}^t_{\BS x},
\end{equation}
\begin{equation}
\left\langle\parderiv{\BS F}{\BS x}\right\rangle(\BS x)=
\int_{\mathcal A_{\BS x}}\parderiv{}{\BS x}\BS F(\BS x,\BS z)
\dif\mu_{\BS x}(\BS z),
\end{equation}
\begin{equation}
\label{eq:H_corr}
\BS H(\BS x)=\int_0^\infty\left(\int_{\mathcal A_{\BS x}}\parderiv{\BS F}
{\BS y}(\BS x,\phi_{\BS x}^s\BS z)\BS T_{\BS x,\BS z}^s\parderiv{\BS G}
{\BS x}(\BS x, \BS z)\dif\mu_{\BS x}(\BS z)\right)\dif s.
\end{equation}
\end{subequations}
Under the ergodicity assumption, one can replace the measure averages
with the time averages:
\begin{subequations}
\begin{equation}
\left\langle\parderiv{\BS F}{\BS x}\right\rangle(\BS x)=
\lim_{r\to\infty}\frac 1r\int_0^r\parderiv{}{\BS x}\BS F(\BS x,
\BS z(\tau))\dif\tau,
\end{equation}
\begin{equation}
\BS H(\BS x)=\int_0^\infty\left(\lim_{r\to\infty}\frac 1r
\int_0^r\parderiv{\BS F}{\BS y}(\BS x,\BS z(\tau+s))\BS T_{\BS x,\BS
  z(\tau)}^s \parderiv{\BS G}{\BS x}(\BS x,\BS
z(\tau))\dif\tau\right)\dif s.
\end{equation}
\end{subequations}

\section{Deterministic linear coupling with total energy conservation}
\label{sec:coupling}

It is often the case that the deterministic parts of the coupling in
$\BS F(\BS x,\BS y)$ and $\BS G(\BS x,\BS y)$ between the slow
variables $\BS x$ and the fast variables $\BS y$ allow the total
energy conservation, as long as in the absence of a random noise the
uncoupled dynamics preserve the energy separately for $\BS x$ and $\BS
y$. As an example, one can look at the model of mean flow -- small
scale interactions via topographic stress in \cite{GroMajGro}, where
the total energy conservation in the coupling between the zonal mean
flow and small scale fluctuations is a key requirement, and topography
can be viewed as a small parameter. Here we assume that the slow and
fast variables are coupled through a deterministic coupling with the
total energy conservation property. More precisely, we let the
functions $\BS F(\BS x,\BS y)$ and $\BS G(\BS x,\BS y)$ take the form
\begin{equation}
\BS F(\BS x,\BS y)=\BS f(\BS x) +\BS f^\prime(\BS x,\BS y),\qquad
\BS G(\BS x,\BS y)=\BS g(\BS y) +\BS g^\prime(\BS x,\BS y),
\end{equation}
and assume that that there exists a symmetric positive definite
quadratic form
\begin{equation}
\label{eq:energy}
E=E_x+\varepsilon E_y,\qquad E_x=\frac 12\BS x^T\BS S_x\BS x, \qquad
E_y=\frac 12\BS y^T\BS S_y\BS y,
\end{equation}
where $\BS S_x$ and $\BS S_y$ are constant symmetric positive definite
matrices, which is preserved by $\BS f^\prime(\BS x,\BS y)$ and $\BS
g^\prime(\BS x,\BS y)$ for any $\BS x$ and $\BS y$:
\begin{equation}
\label{eq:en_cons}
\BS x^T\BS S_x\BS f^\prime(\BS x,\BS y)+\BS y^T\BS S_y\BS g^\prime(\BS
x,\BS y)=0.
\end{equation}
For \eqref{eq:en_cons} to hold, it is sufficient to have $\BS
f^\prime(\BS x,\BS y)$ and $\BS g^\prime(\BS x,\BS y)$ in the form
\begin{subequations}
\label{eq:fp_gp}
\begin{equation}
\BS f^\prime(\BS x,\BS y)=\BS S_x^{-1/2}\BS L(\BS x,\BS y)\BS
S_y^{1/2}\BS y,
\end{equation}
\begin{equation}
\BS g^\prime(\BS x,\BS y)=-\BS S_y^{-1/2}\BS L^T(\BS x,\BS y)\BS
S_x^{1/2}\BS x,
\end{equation}
\end{subequations}
for an arbitrary $N_x\times N_y$ matrix-valued function $\BS L(\BS
x,\BS y)$. Additionally, if $\BS f^\prime(\BS x,\BS y)$ and $\BS
g^\prime(\BS x,\BS y)$ are analytic on $\mathbb R^{N_x}\times\mathbb
R^{N_y}$, the relations in \eqref{eq:fp_gp} are also necessary:
indeed, in this case $\BS f^\prime(\BS x,\BS y)$ and $\BS g^\prime(\BS
x,\BS y)$ can be represented as
\begin{equation}
\BS f^\prime(\BS x,\BS y)=\BS f^\prime(\BS x,\BS 0)+\BS f_1^\prime(\BS
x,\BS y)\BS y,\qquad \BS g^\prime(\BS x,\BS y)=\BS g^\prime(\BS 0,\BS
y)+\BS g_1^\prime(\BS x,\BS y)\BS x,
\end{equation}
where $\BS f_1^\prime(\BS x,\BS y)\BS y$ and $\BS g_1^\prime(\BS x,\BS
y)\BS x$ represent the Taylor expansion terms above zero order in $\BS
y$ and $\BS x$, respectively. As a result, \eqref{eq:en_cons} becomes
\begin{equation}
\BS x^T\left[\BS S_x\BS f_1^\prime(\BS x,\BS y)+\BS
g_1^\prime(\BS x,\BS y)^T\BS S_y\right]\BS y=-\BS x^T\BS S_x\BS f^\prime
(\BS x,\BS 0)-\BS y^T\BS S_y\BS g^\prime(\BS 0,\BS y),
\end{equation}
where setting either $\BS x=\BS 0$ or $\BS y=\BS 0$ turns the
left-hand side to zero, which means that $\BS f^\prime(\BS x,\BS
0)=\BS g^\prime(\BS 0,\BS y)=\BS 0$ for any $\BS x$ and $\BS y$, and
that $\BS f^\prime(\BS x,\BS y)=\BS f_1^\prime(\BS x,\BS y)\BS y$,
$\BS g^\prime(\BS x,\BS y)=\BS g_1^\prime(\BS x,\BS y)\BS x$. This
results in
\begin{equation}
\BS S_x\BS f_1^\prime(\BS x,\BS y)=-\BS g_1^\prime(\BS x,\BS y)^T\BS S_y.
\end{equation}
Let $\BS L(\BS x,\BS y)$ and $\BS L^*(\BS x,\BS y)$ be given,
respectively, by
\begin{equation}
\BS L(\BS x,\BS y)=\BS S_x^{1/2}\BS f_1^\prime(\BS x,\BS y)\BS
S_y^{-1/2},\qquad\BS L^*(\BS x,\BS y)=\BS S_y^{1/2}\BS g_1^\prime (\BS
x,\BS y)\BS S_x^{-1/2},
\end{equation}
then, obviously, $\BS L^*(\BS x,\BS y)=-\BS L(\BS x,\BS y)^T$ for any
$\BS x$ and $\BS y$. Therefore, we end up with
\begin{equation}
\BS f_1^\prime(\BS x,\BS y)=\BS S_x^{-1/2}\BS L(\BS x,\BS y)\BS
S_y^{1/2},\qquad \BS g_1^\prime(\BS x,\BS y)=-\BS S_y^{-1/2}\BS
L^T(\BS x,\BS y)\BS S_x^{1/2}
\end{equation}
for some $N_x\times N_y$ matrix-valued function $\BS L(\BS x,\BS y)$,
and \eqref{eq:fp_gp} follows directly.

In this work we study the simplest case of the energy-preserving
coupling with $\BS L(\BS x,\BS y)=\BS L$, where $\BS L$ is a constant
matrix and does not depend on the variables $\BS x$ and $\BS y$, so
that the coupling is linear with respect to $\BS y$ (for the slow
variables) and $\BS x$ (for the fast variables). With
\eqref{eq:fp_gp}, \eqref{eq:dyn_sys} becomes
\begin{subequations}
\label{eq:dyn_sys_en}
\begin{equation}
\label{eq:dyn_sys_slow_en}
\deriv{\BS x}t=\BS f(\BS x)+\BS S_x^{-1/2} \BS L\BS S_y^{1/2}\BS
y,
\end{equation}
\begin{equation}
\label{eq:dyn_sys_fast_en}
\deriv{\BS y}t=\frac 1\varepsilon\left[\BS g(\BS y)-\BS S_y^{-1/2}
\BS L^T\BS S_x^{1/2}\BS x\right],
\end{equation}
\end{subequations}
with the limiting fast dynamics for constant $\BS x$ given by
\begin{equation}
\label{eq:dyn_sys_fast_limiting_z_en}
\deriv{\BS z}t=\BS g(\BS z)-\BS S_y^{-1/2}\BS L^T\BS S_x^{1/2}\BS x.
\end{equation}
With the above changes, the integral of the time autocorrelation
function in \eqref{eq:H_corr} can be written as
\begin{equation}
\label{eq:H}
\begin{split}
\BS H(\BS x)=-\BS S_x^{-1/2}\BS L\BS S_y^{1/2}
\left(\int_0^\infty\int_{\mathcal A_{\BS x}}\BS T_{\BS x,\BS
  z}^s\dif\mu_{\BS x}(\BS z)\dif s\right) \BS S_y^{-1/2}\BS L^T\BS
S_x^{1/2}.
\end{split}
\end{equation}
The dynamics for the averaged slow tangent map from
\eqref{eq:dyn_sys_slow_tangent} is now written as
\begin{equation}
\label{eq:dyn_sys_slow_tangent_2}
\begin{split}
\parderiv{}t\BS{TX}^t_{\BS x}=\left(\parderiv{}{\BS x}\BS f(\BS
x(t))+\BS H(\BS x(t))\right)\BS{TX}^t_{\BS x}.
\end{split}
\end{equation}
Observe that if the slow system in \eqref{eq:dyn_sys_slow_en} is
completely decoupled from the fast variables (for example, $\BS L=0$),
then the linear tangent model above loses the term $\BS H(\BS x)$,
otherwise remaining the same. Below, we study the effect of $\BS H(\BS
x)$ on the dynamics of \eqref{eq:dyn_sys_slow_tangent_2}. Observe that
$\BS H(\BS x)$ can be written in the form
\begin{subequations}
\label{eq:H_C_C}
\begin{equation}
\BS H(\BS x)=-\BS S_x^{-1/2}\BS L\BS S_y^{1/2}\BS{\mathcal C}(\BS
x)\BS S_y^{-1/2}\BS L^T\BS S_x^{1/2},
\end{equation}
\begin{equation}
\label{eq:corr_1}
\BS{\mathcal C}(\BS x)=\int_0^\infty\BS C(\BS x,s)\dif s,
\end{equation}
\begin{equation}
\label{eq:corr_2}
\BS C(\BS x,s)=\int_{\mathcal A_{\BS x}} \BS T_{\BS x,\BS z}^s
\dif\mu_{\BS x}(\BS z)=\lim_{r\to\infty}\frac 1r\int_0^r\BS T_{\BS x,
\BS z(\tau)}^s\dif\tau
\end{equation}
\end{subequations}
(the last equality is obtained under the ergodicity assumption), where
$\BS{\mathcal C}(\BS x)$ is the infinite time linear response operator
of the mean state of \eqref{eq:dyn_sys_fast_limiting_z_en} to a
constant external forcing, with given parameter $\BS x$ (for details,
see Appendix \ref{sec:inf_time_add}). Replacing $\BS{\hat x}=\BS
S_x^{1/2}\BS x$ and $\BS{\hat y}=\BS S_y^{1/2}\BS y$, which are
canonical variables for the energy in \eqref{eq:energy} (that is, the
quadratic form in \eqref{eq:energy} is diagonal in the new variables
$\BS{\hat x}$ and $\BS{\hat y}$), in the new variables we obtain
\begin{equation}
\BS{\hat H}(\BS x)=-\BS L\BS{\hat{\mathcal C}}(\BS x)\BS L^T.
\end{equation}
Here,
\begin{equation}
\BS{\hat H}(\BS x)=\BS S_x^{1/2}\BS H(\BS x)\BS S_x^{-1/2}
\end{equation}
is the term $\BS H(\BS x)$ in the canonical energy coordinates, while
$\BS{\hat{\mathcal C}}(\BS x)$ is the infinite-time mean state linear
response operator for \eqref{eq:dyn_sys_fast_limiting_z_en} in the
canonical energy coordinates $\BS{\hat x}$ and $\BS{\hat y}$:
\begin{equation}
\BS{\hat{\mathcal C}}(\BS x)=\BS S_y^{1/2}\BS{\mathcal C}(\BS x)\BS
S_y^{-1/2},
\end{equation}
since for $\BS{\hat z}=\BS S_y^{1/2}\BS z$ we have
\begin{equation}
\BS{\hat T}_{\BS x,\BS{\hat z}}^s=\BS S_y^{1/2}
\BS T_{\BS x,\BS z}^s\BS S_y^{-1/2}.
\end{equation}

\section{Suppression of chaos at slow variables by the linear
energy-preserving coupling for the rapidly decorrelating fast
dynamics}
\label{sec:suppression}

Without loss of generality, further we assume that $\BS x$ and $\BS y$
are already the canonical energy variables $\BS{\hat x}$ and $\BS{\hat
  y}$, as if the dynamical system in \eqref{eq:dyn_sys_en} was
formulated in the canonical energy variables from the beginning, such
that $\BS S_x$ and $\BS S_y$ in \eqref{eq:energy} are multiples of an
identity matrix. With this assumption, the term $\BS H(\BS x)$ from
\eqref{eq:H} is written as
\begin{equation}
\label{eq:H2}
\BS H(\BS x)=-\BS L\BS{\mathcal C}(\BS x)\BS L^T,
\end{equation}
where $\BS{\mathcal C}(\BS x)$ is given by \eqref{eq:corr_1} and
\eqref{eq:corr_2}, computed from the fast limiting dynamical system
given by
\begin{equation}
\label{eq:dyn_sys_fast_limiting_z_en_canon}
\deriv{\BS z}t=\BS g(\BS z)-\BS L^T\BS x,
\end{equation}
with given parameter $\BS x$.

At this point, observe that if $\BS H$ is negative definite, then it
has a damping effect on the growth of $\BS{TX}^t_{\BS x}$ in
\eqref{eq:dyn_sys_slow_tangent_2} (more generally, in the
non-canonical variables $\BS H$ reduces the slow energy $E_x$,
promoting the Lyapunov stability). However, from \eqref{eq:H2} we see
that $\BS H$ can only be negative definite whenever $\BS{\mathcal C}$
is positive definite. Indeed, let $\BS v^T\BS H\BS v<0$ for any
nonzero $\BS v$, then, denoting $\BS w=\BS L^T\BS v$, we obtain $\BS
w^T\BS{\mathcal C}\BS w>0$. Below we are going to establish the
physical conditions which promote the positive-definiteness of
$\BS{\mathcal C}(\BS x)$.

First, observe that $\BS{\mathcal C}(\BS x)$ is the infinite-time
linear response operator for
\eqref{eq:dyn_sys_fast_limiting_z_en_canon} perturbed by a small
constant forcing. In particular, if
\eqref{eq:dyn_sys_fast_limiting_z_en_canon} is an Ornstein-Uhlenbeck
process \cite{OrnUhl} of the form
\begin{equation}
\dif\BS z=\left(\BS h-\BS\Gamma\BS z-\BS L^T\BS x\right)\dif
s+\BS\sigma\dif\BS W_s,
\end{equation}
with $\BS\Gamma$ being positive-definite matrix to (almost) ensure
boundedness of the solution, then we obtain explicitly
\begin{subequations}
\begin{equation}
\BS C(\BS x,s)=\BS C(s)=e^{-s\BS\Gamma},
\end{equation}
\begin{equation}
\BS{\mathcal C}(\BS x)=\BS{\mathcal C}=\int_0^\infty \BS
C(s)\dif s=\BS\Gamma^{-1},
\end{equation}
\begin{equation}
\BS v^T\BS H\BS v=-\BS v^T\BS L\BS\Gamma^{-1}\BS L^T\BS v<0,
\end{equation}
\end{subequations}
which means that the term $\BS H$ automatically contributes towards
the reduction of chaos at the slow variables whenever the fast
dynamics are modeled by an Ornstein-Uhlenbeck process. In general,
since
\begin{equation}
\langle\delta\BS z\rangle=\BS{\mathcal C}\cdot\delta\BS f,
\end{equation}
where $\delta\BS f$ is a small constant perturbation, and
$\langle\delta\BS z\rangle$ is the infinite-time linear response of
the mean state $\langle\BS z\rangle$ of
\eqref{eq:dyn_sys_fast_limiting_z_en_canon}, the positive-definiteness
of $\BS{\mathcal C}$ means that
\begin{equation}
\langle\delta\BS z\rangle\cdot\delta\BS f>0
\mbox{ for all sufficiently small }\delta\BS f,
\end{equation}
that is, the response of the mean state never develops against the
perturbation. It is not difficult to show that the following identity
holds whenever $\delta\BS f$ vanishes:
\begin{equation}
\deriv{}t\langle\delta\BS z\rangle=-\delta\BS f,
\end{equation}
that is, the time derivative of $\langle\delta\BS z\rangle$ at the
moment the small constant perturbation $\delta\BS f$ vanishes equals
this perturbation with the opposite sign. Therefore, the
positive-definiteness of $\BS{\mathcal C}$ is equivalent to the
following mean stability property of
\eqref{eq:dyn_sys_fast_limiting_z_en_canon}:
\begin{equation}
\label{eq:mean_stab}
\deriv{}t\|\langle\delta\BS z\rangle\|<0,
\end{equation}
that is, any sufficiently small infinite-time perturbation of the mean
state $\langle\delta\BS z\rangle$ decreases in time at the moment when
the external perturbation $\delta\BS f$ is removed. At present, it is
not precisely clear to the author how to ascertain this property for
general nonlinear $\BS g$ in
\eqref{eq:dyn_sys_fast_limiting_z_en_canon}, however, one might expect
such properties to be common in strongly chaotic and mixing turbulent
dynamics, as their statistical properties are often modeled by an
appropriate Ornstein-Uhlenbeck process (see
\cite{MajTimVan,MajTimVan2,MajTimVan3} and references therein), and,
as shown above, the Ornstein-Uhlenbeck process satisfies this property
automatically. For general nonlinear dynamics, the mean stability
property must be associated with the situation where the typical
Poincar\'e recurrence time of nonlinear motion around the mean state
in \eqref{eq:dyn_sys_fast_limiting_z_en_canon} (which can be viewed as
an advective time scale) is not much shorter than the turbulent mixing
autocorrelation time. The reason is that, since $\BS T^0_{\BS z}=\BS
I$ (the identity matrix) for any $\BS z$, then there always exists
$a^*$ such that
\begin{equation}
\BS{\mathcal C}_a(\BS x)=\int_0^{a^*}\BS C(\BS x,s)\dif s
\mbox{ is positive definite for all }a,\quad 0<a<a^*.
\end{equation}
The positive-definiteness of $\BS{\mathcal C}_a(\BS x)$ for larger $a$
can be violated by the domination of the rotation part in $\BS C(\BS
x,s)$ for larger $s$, which evolves on the advective time scale of
\eqref{eq:dyn_sys_fast_limiting_z_en_canon}. However, this effect can
be prevented by a rapid decay of $\|\BS C(\BS x,s)\|$ for large $s$,
which is governed by the turbulent mixing autocorrelation time. Thus,
in general, one can expect the positive-definiteness of $\BS{\mathcal
  C}(\BS x)$ in the situations where the turbulent mixing
autocorrelation time scale is not much longer than the advective time
scale.

\section{Revisiting the rescaled Lorenz model}
\label{sec:lorenz_revisited}

Here, after developing the theory for the suppression of chaos at the
slow variables by the rapidly mixing fast dynamics, we return back to
the rescaled Lorenz model in \eqref{eq:lorenz_rescaled}. At this
point, we can observe that the rescaled Lorenz model in
\eqref{eq:lorenz_rescaled} conforms to the requirements of the theory
developed above, with the following properties:
\begin{itemize}
\item The $\BS x$ and $\BS y$ variables in the rescaled Lorenz model
  are already the canonical energy variables (up to constant positive
  factors). Indeed, set all constant and linear terms in the rescaled
  Lorenz model in \eqref{eq:lorenz_rescaled} to zero. Then, it is easy
  to see that the nonlinear parts in $\BS x$ and $\BS y$ separately
  preserve the energies
  \begin{equation}
    E_x=\frac{\lambda_x}2\sum_ix_i^2,\qquad E_y=\frac{\lambda_y}{2J}
    \sum_{i=1}^{N_x}\sum_{j=1}^Jy_{ij}^2,
  \end{equation}
  while the full rescaled Lorenz model without forcing and
  dissipation preserves the total energy of the form
  \begin{equation}
    E=\frac{\lambda_x}2\sum_{i=1}^{N_x}x_i^2+\frac{\varepsilon\lambda_y}
    {2J}\sum_{i=1}^{N_x}\sum_{j=1}^Jy_{ij}^2=E_x+\varepsilon E_y,
  \end{equation}
  Note that $\BS x$ and $\BS y$ in the rescaled Lorenz model are
  already the canonical energy variables, as the quadratic form above
  is diagonal.
\item The deterministic linear energy-preserving coupling is given by
  \begin{equation}
    (\BS L\BS y)_i=-\sum_{j=1}^Jy_{ij},\qquad (\BS L^T\BS
    x)_{ij}=-x_i.
  \end{equation}
  Then, $\BS H(\BS x)$ is given by
  \begin{equation}
    \BS H(\BS x)=-\frac{\lambda_x\lambda_y}J\BS L\BS{\mathcal C}(\BS
    x)\BS L^T.
  \end{equation}
\item The rescaled Lorenz model can be written as
  \begin{equation}
    \deriv{\BS x}t=\BS f(\BS x)+\frac{\lambda_y}J\BS L\BS y,
  \end{equation}
  \begin{equation}
    \deriv{\BS y}t=\frac 1\varepsilon\left[\BS g(\BS y)- \lambda_x\BS
      L^T\BS x\right],
  \end{equation}
  with $\BS f(\BS x)$ and $\BS g(\BS y)$ given by
  \begin{subequations}
    \begin{equation}
      f_i(\BS x)=x_{i-1}(x_{i+1}-x_{i-2})+\frac 1{\beta_x} \left(\bar
      x(x_{i+1}-x_{i-2})-x_i\right)+\frac{F_x-\bar x}{\beta_x^2},
    \end{equation}
    \begin{equation}
      g_{ij}(\BS y)=y_{i,j+1}(y_{i,j-1}-y_{i,j+2})+\frac 1{\beta_y}
      \left(\bar y(y_{i,j-1}-y_{i,j+2})-y_{i,j}\right)+\frac{F_y-\bar
        y} {\beta_y^2},
    \end{equation}
  \end{subequations}
  respectively.
\item The slow tangent model dynamics for the rescaled Lorenz model
  can be written as
  \begin{equation}
  \deriv{}t\BS{TX}^t_{\BS x}=\left[\parderiv{\BS f}{\BS x}(\BS
    x(t))+\BS H(\BS x(t))\right]\BS{TX}^t_{\BS x}.
  \end{equation}
\end{itemize}
Overall, it turns out that in the slow-time tangent linear model for
the rescaled Lorenz system in \eqref{eq:lorenz_rescaled}, the
additional term $\BS H(\BS x)$ suppresses chaos at the slow variables
if the infinite-time linear response operator $\BS{\mathcal C}(\BS x)$
of the limiting fast dynamics, given by \eqref{eq:corr_1} and
\eqref{eq:corr_2}, is positive definite. It is, of course, not
practically feasible to compute $\BS C(\BS x,s)$ directly for all $\BS
x$ and $s$. Here, however, we will use a suitable approximation for
$\BS C(\BS x,s)$ in the form of the quasi-Gaussian linear response
operator
\cite{Abr5,Abr6,Abr7,AbrMaj4,AbrMaj5,AbrMaj6,AbrMaj7,MajAbrGro}. Namely,
we approximate the average above by the average with respect to the
Gaussian distribution $p_{\BS x}^G(\BS z)$ with the same mean state
and covariance matrix as $\mu_{\BS x}$ \cite{AbrMaj4,AbrMaj5,AbrMaj6}:
\begin{subequations}
\begin{equation}
\label{eq:C_lorenz_2}
\BS C(\BS x,s)\approx\BS C_G(\BS x,s)=\int_{\mathbb R^{N_y}}\BS
T_{\BS x,\BS z}^sp_{\BS x}^G(\BS z)\dif\BS z,
\end{equation}
\begin{equation}
p_{\BS x}^G(\BS z)=\exp\left(-\frac 12(\BS z-\BS{\bar z}(\BS x))^T
\BS\Sigma^{-1}(\BS x)(\BS z-\BS{\bar z}(\BS x))\right),
\end{equation}
\end{subequations}
where $\BS{\bar z}(\BS x)$ and $\BS\Sigma(\BS x)$ are the mean state
and covariance matrix, respectively, of the fast limiting dynamics in
\eqref{eq:dyn_sys_fast_limiting_z_en_canon} for fixed $\BS x$. Then,
via integration by parts, one can rewrite \eqref{eq:C_lorenz_2} as
\begin{equation}
\BS C_G(\BS x,s)=\left[\int_{\mathbb R^{N_y}}\BS z(s)(\BS z-
\BS{\bar z}(\BS x))^Tp_{\BS x}^G(\BS z)\dif\BS z\right]
\BS\Sigma^{-1}(\BS x),
\end{equation} 
which after going back to the time averages becomes the time
autocorrelation function of the form
\begin{equation}
\BS C_G(\BS x,s)=\left[\lim_{r\to\infty}\frac 1r\int_0^r\BS z(\tau+s)
(\BS z(\tau)-\BS{\bar z}(\BS x))^T\dif\tau\right]\BS\Sigma^{-1}(\BS x),
\end{equation}
which, after denoting $\BS z(t)=\BS{\bar z}(\BS x)+\BS z^\prime(t)$, can be
written in the form
\begin{equation}
\label{eq:C_lorenz_3}
\BS C_G(\BS x,s)=\left[\lim_{r\to\infty}\frac 1r\int_0^r\BS z^\prime(\tau+s)
\BS z^{\prime T}(\tau)\dif\tau\right]\BS\Sigma^{-1}(\BS x).
\end{equation}
Even after the above simplification, it is still not practically
feasible to compute $\BS C(\BS x,s)$ at many points $\BS x$, as needed
in \eqref{eq:dyn_sys_slow_tangent_2}. Here, however, we only need to
observe $\BS C(\BS x,s)$ and $\BS{\mathcal C}(\BS x)$ for diagnostic
purposes, to relate to the trends in Figures \ref{fig:av_diff_1} and
\ref{fig:av_diff_2}. Thus, we compute $\BS{\mathcal C}_G(\BS x)$ at
the single point $\BS x=\BS{\bar x}$, which is the long-term mean
state of the slow variables in \eqref{eq:lorenz_rescaled}. Since the
motion of the slow variables occurs in the vicinity of their mean
state, the computed $\BS{\mathcal C}_G(\BS{\bar x})$ should generally
reflect the trends happening around that point.
\begin{figure}%
\picturehere{\figdir/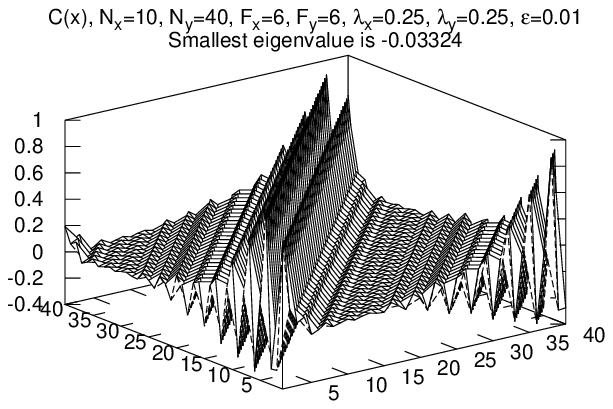}%
\picturehere{\figdir/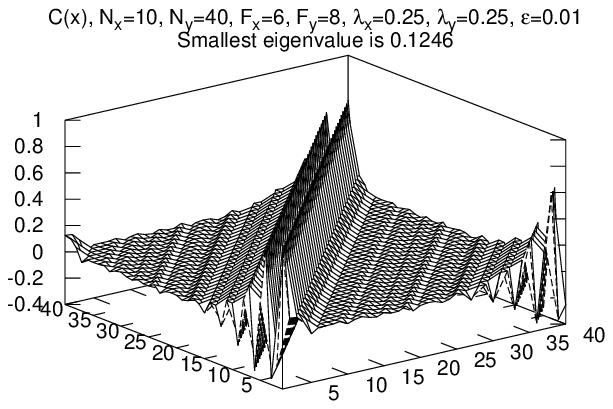}\\%
\picturehere{\figdir/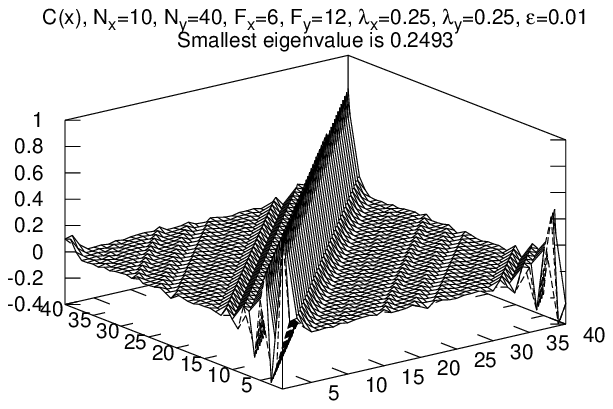}%
\picturehere{\figdir/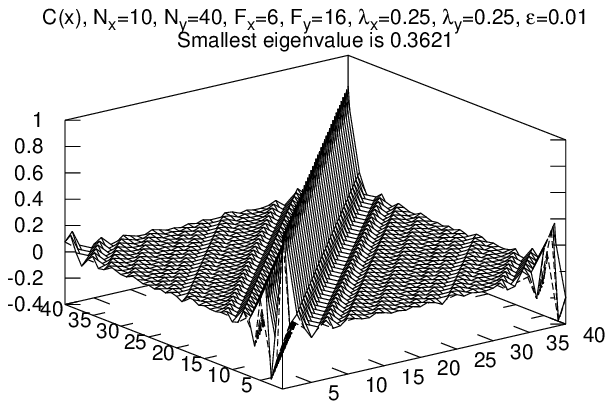}\\%
\begin{center}%
\picturehere{\figdir/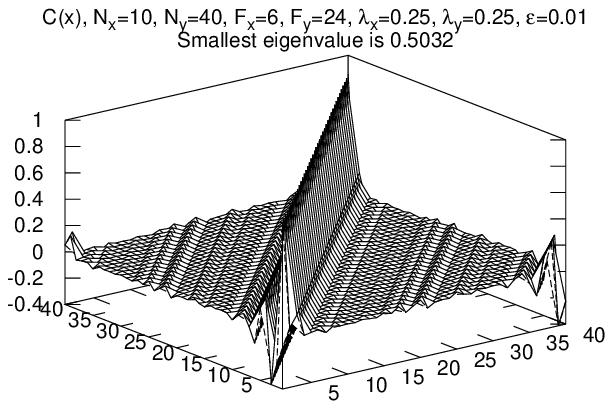}%
\end{center}%
\caption{Quasi-Gaussian approximations of averaged infinite-time
  linear response operators $\BS{\mathcal C}_G(\BS{\bar x})$ for the
  rescaled Lorenz model in \eqref{eq:lorenz_rescaled} with the
  following parameters: $N_x=10$, $N_y=40$, $F_x=6$, $F_y=6,8,12,16$
  and $24$, $\lambda_x=\lambda_y=0.25$, $\varepsilon=0.01$.}
\label{fig:inf_time_resp_1}
\end{figure}%
In Figure \ref{fig:inf_time_resp_1} we show the computed operators
$\BS{\mathcal C}_G(\BS{\bar x})$, together with the smallest
eigenvalues of the symmetric parts of those operators (as
skew-symmetric parts do not contribute to positive-definiteness).
Obviously, the larger is the smallest eigenvalue of a symmetric part
of the matrix, the ``more positive-definite'' is the matrix
itself. Indeed, observe that the smallest eigenvalue of the response
systematically increases with increased $F_y$ (in fact, for $F_y=6$
the operator is not even positive definite). Of course, these results
are no more than crude estimates of the actual trends, as in reality
the response operators must be computed for each $\BS x$ and using the
exact formulas in \eqref{eq:corr_1} and \eqref{eq:corr_2}, rather than
the one-point estimate under Gaussian assumption. Yet, as one can see,
even these crude approximations help to connect the trends in Figure
\ref{fig:inf_time_resp_1} with the surprising behavior displayed in
Figures \ref{fig:av_diff_1} and \ref{fig:av_diff_2} in Section
\ref{sec:lorenz}.
\begin{figure}%
\begin{center}%
\picturehere{\figdir/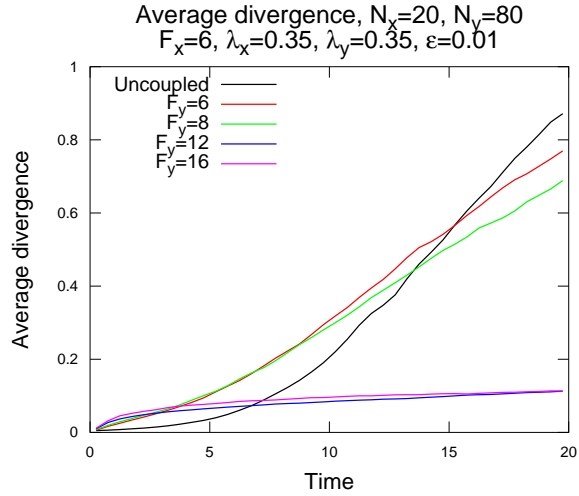}%
\end{center}%
\caption{Average divergence between perturbed and unperturbed running
  averages of the slow variables for the rescaled Lorenz model in
  \eqref{eq:lorenz_rescaled} with the following parameters: $N_x=20$,
  $N_y=80$, $F_x=6$, $F_y=6,8,12$ and $16$,
  $\lambda_x=\lambda_y=0.35$, $\varepsilon=0.01$, as well as for the
  uncoupled rescaled Lorenz model with $N=20$ and $F=6$.}%
\label{fig:av_diff_3}%
\end{figure}%
\begin{figure}%
\largepicturehere{\figdir/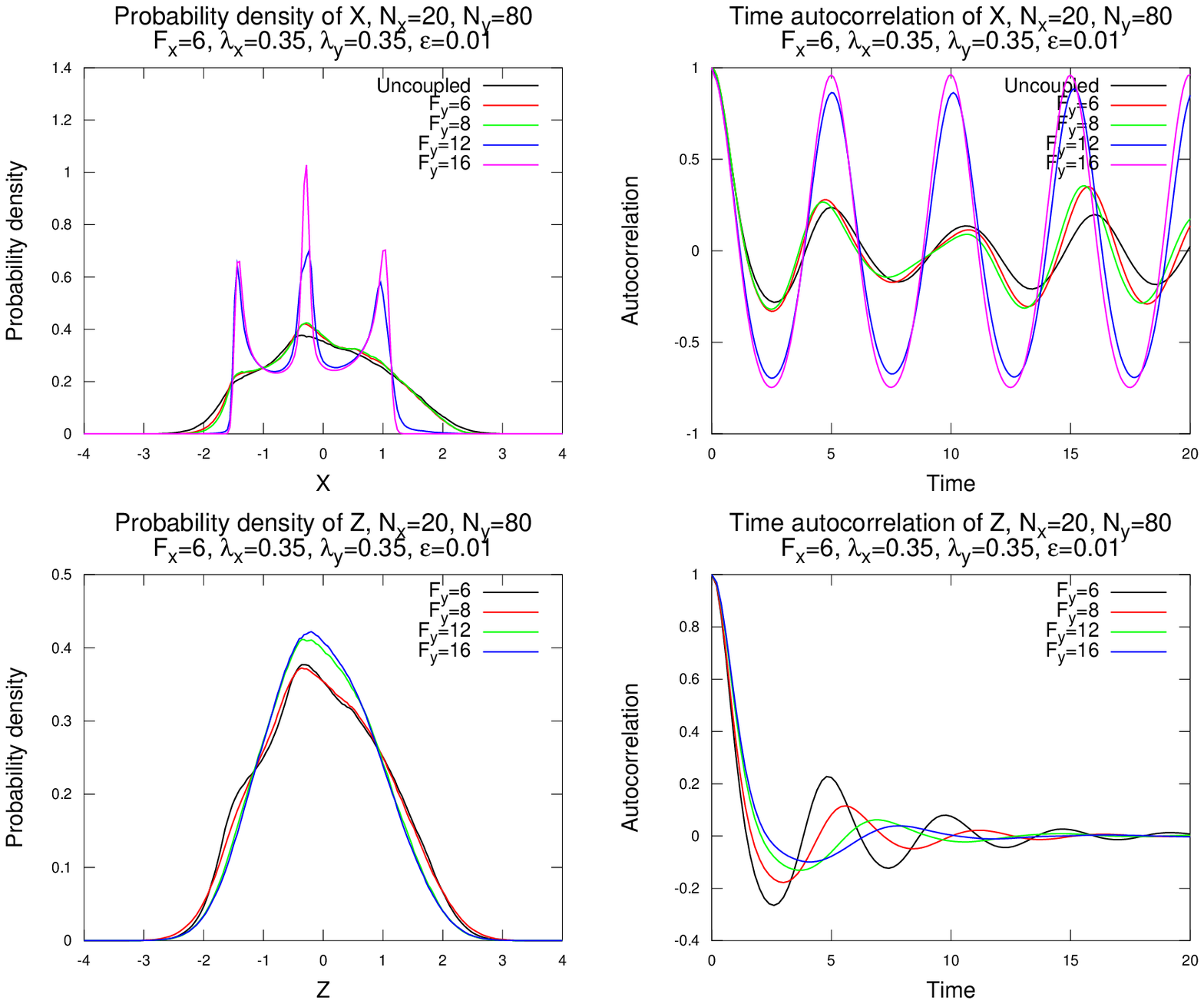}%
\caption{Probability density functions and time autocorrelation
  functions for $\BS x$ of \eqref{eq:lorenz_rescaled} and $\BS z$ of
  \eqref{eq:lorenz_fast_limiting} (with free parameter $\BS x$ set to
  its mean value) with the following parameters: $N_x=20$, $N_y=80$,
  $F_x=6$, $F_y=6,8,12$ and $16$, $\lambda_x=\lambda_y=0.35$,
  $\varepsilon=0.01$.}%
\label{fig:pdf_corr_2}%
\end{figure}%
\begin{table}%
\begin{center}%
\begin{tabular}{|c|}%
\hline
$N_x=20$, $N_y=80$, $F_x=6$, $\lambda_x=\lambda_y=0.35$, $\varepsilon=0.01$ \\
\hline
\begin{tabular}{c||c|c|c|c}
$F_y$ & $\BS x$-mean & $\BS x$-var & $\BS y$-mean & $\BS y$-var \\
\hline\hline
$6$ & $6.216\cdot 10^{-3}$ & $0.8878$ & $-9.982\cdot 10^{-3}$ & $1.119$ \\
$8$ & $2.34\cdot 10^{-2}$ & $0.8728$ & $-2.791\cdot 10^{-2}$ & $1.19$ \\
$12$ & $-0.1363$ & $0.6927$ & $-5.553\cdot 10^{-2}$ & $1.168$ \\
$16$ & $-0.1444$ & $0.6703$ & $-8.669\cdot 10^{-2}$ & $1.199$ \\
\end{tabular}\\
\hline
\end{tabular}%
\end{center}%
\caption{The mean states and variances of the $\BS x$ and $\BS y$
  variables for the rescaled Lorenz model in
  \eqref{eq:lorenz_rescaled} with the following parameters: $N_x=20$,
  $N_y=80$, $F_x=6$, $F_y=6,8,12$ and $16$,
  $\lambda_x=\lambda_y=0.35$, $\varepsilon=0.01$.}
\label{tab:mean_var_2}
\end{table}

\section{Complete suppression of chaos at slow variables while
the uncoupled slow dynamics remain chaotic}
\label{sec:suppression_uncoupled}

One of the key questions in the atmosphere/ocean science is whether
the uncoupled system, consisting of slow variables only, is more or
less chaotic than its original version, coupled with the fast, often
unresolved or underresolved variables. Indeed, often scientists work
with uncoupled models consisting of slow variables only, where
coupling terms were replaced with the estimates of the long-term
averages of the corresponding fast variables (for example, the T21
barotropic model with the realistic Earth topography
\cite{AbrMaj6,Fra,Sel}), and study dynamical properties of the slow
variables in the uncoupled models. The common sense in this case
suggests that if the uncoupled slow model is chaotic, then, naturally,
its original version coupled with fast rapidly mixing dynamics should
be even more chaotic.

Remarkably, the common sense logic in this situation is deceiving. In
fact, it turns out to be possible even to reach the transition from
the chaotic to stable slow dynamics by increasing the turbulent mixing
at the fast variables, while the uncoupled slow dynamics remain
chaotic. This generally follows from the fact that in
\eqref{eq:dyn_sys_slow_tangent_2} the term $\BS H$ is absent from the
uncoupled dynamics, and only the Jacobian of $\BS f$ remains. Since
$\BS H$ acts as a dissipation in \eqref{eq:dyn_sys_slow_tangent_2} for
positive-definite infinite-time linear response of the fast variables,
it is quite natural that in its absence the uncoupled dynamics is more
chaotic than the fully coupled dynamics at the slow variables.

Here we demonstrate such an example for the rescaled Lorenz model in
\eqref{eq:lorenz_rescaled} with the following set of parameters:
$N_x=20$, $N_y=80$, $F_x=6$, $\lambda_x=\lambda_y=0.35$,
$\varepsilon=0.01$, and compare it with the uncoupled rescaled Lorenz
model in \eqref{eq:L96_rescaled} with $N=20$ and $F=6$. In Figure
\ref{fig:av_diff_3} we show the average divergence of a perturbed
trajectory from an unperturbed one for this set of parameters. Observe
that, while the slow dynamics for $F_y=6,8$ and for uncoupled dynamics
with $F=6$ are clearly chaotic, for $F_y=12$ and greater values the
abrupt transition occurs, where the difference between the perturbed
and unperturbed $\BS x$ time series does not grow much beyond
10\%. The statistics for the same set of parameters are shown in
Figure \ref{fig:pdf_corr_2} and Table \ref{tab:mean_var_2}.  Here,
clearly the correlations and probability densities for the slow
variables differ significantly for the chaotic and stable regimes,
with the statistics for the uncoupled Lorenz model with $F=6$ being
almost identical to those for the chaotic regimes with
$F_y=6,8$. However, the chaos at slow variables is suppressed purely
by the dynamical mechanism uncovered in this work; indeed, the mean
state and variance of both the slow and fast variables do not change
substantially enough to suppress chaos by creating a counteracting
mean forcing term at the slow dynamics to suppress $F_x$ (see Table
\ref{tab:mean_var_2} for confirmation), while the time autocorrelation
functions for the fast variables in Figure \ref{fig:pdf_corr_2} with
$\BS x$ set to the mean state decay faster for larger values of $F_y$,
indicating stronger mixing. The key observation here is that the
behavior of the uncoupled model with only slow variables in the same
regime is deceiving -- it is clearly chaotic, while the full two-scale
model loses chaos at the slow variables in more turbulent regimes of
the fast dynamics.

\section{Generalization onto the two-scale setting with stochastic
noise and explicit time dependence}
\label{sec:extension}

For the simplicity of the presentation, the theory above was developed
for a simplified setting in \eqref{eq:dyn_sys}, which matched the
structure of the two-scale rescaled Lorenz model. However, the
real-world chaotic multiscale models often include the stochastic
noise terms. Additionally, the explicit dependence on the range of
different time scales is quite common in practice. As an example, one
can look at the annual and diurnal time scales of the external solar
forcing at a given physical location on the surface of the
Earth. These time scales differ by two orders of magnitude. Below we
develop the framework for a more general setting, which includes both
the stochastic noise and explicit dependence on time.

Here, we modify the dynamical system in \eqref{eq:dyn_sys} by
including the explicit time dependence on both the slow time scale $t$
and the fast time scale $t/\varepsilon$. Additionally, we include the
noise terms at both the slow and fast variables. The dynamical system
in \eqref{eq:dyn_sys} becomes a system of the It\^o stochastic
differential equations of the form
\begin{subequations}
\label{eq:dyn_sys_nonauto}
\begin{equation}
\label{eq:dyn_sys_nonauto_slow}
\dif\BS x=\BS F(\BS x,\BS y,t,t/\varepsilon)\dif t+\BS\kappa(\BS
x,t,t/\varepsilon)\dif\BS V_t,
\end{equation}
\begin{equation}
\label{eq:dyn_sys_nonauto_fast}
\dif\BS y=\frac 1\varepsilon\BS G(\BS x,\BS y,t,t/\varepsilon)\dif t+
\frac 1{\sqrt\varepsilon}\BS\sigma(\BS y,t,t/\varepsilon)\dif\BS W_t.
\end{equation}
\end{subequations}
Above, the terms $\BS V_t$ and $\BS W_t$ are $K_x$- and
$K_y$-dimensional Wiener processes, respectively, for some positive
integers $K_x$ and $K_y$, while $\BS\kappa(\BS x,t,t/\varepsilon)$ and
$\BS\sigma(\BS y,t,t/\varepsilon)$ are the $N_x\times K_x$ and
$N_y\times K_y$ matrix-valued functions, respectively. Observe that
here $\BS\kappa$ does not depend on $\BS y$, while $\BS\sigma$ does
not depend on $\BS x$, as the coupling between the slow and fast
variables must be deterministic to allow the total energy
conservation.

In Section \ref{sec:twoscale} we assumed that both the slow and fast
dynamics in \eqref{eq:dyn_sys} do not depend explicitly on time, and
that assumption rendered the fast limiting dynamics in
\eqref{eq:dyn_sys_fast_limiting_z} into autonomous with the invariant
distribution measure $\mu_{\BS x}$. Since the dynamics in
\eqref{eq:dyn_sys_nonauto} explicitly depend on $t/\varepsilon$, the
fast limiting dynamics as $\varepsilon\to\infty$ (the analogue of that
in \eqref{eq:dyn_sys_fast_limiting_z}) is non-autonomous:
\begin{equation}
\label{eq:dyn_sys_nonauto_fast_limiting_z}
\dif\BS z=\BS G(\BS x,\BS z,t,s)\dif s+\BS\sigma(\BS z,t,s)\dif\BS W_s.
\end{equation}
Here, $\BS x$ and $t$ are the constant parameters, and the invariant
measure $\mu$ for the dynamics of $\BS z(s)$ above may no longer exist
because of the explicit dependence of
\eqref{eq:dyn_sys_nonauto_fast_limiting_z} on $s$. However, practice
shows that, from the physical perspective, averages exist even in a
non-autonomous system (where the averaging is assumed to be with
respect to the time $s$). Usually, the reason for that is that the way
$\BS F$, $\BS G$, $\BS\kappa$ and $\BS\sigma$ depend on $s$ is
bounded, so, in the limit as $s\to\infty$, these functions ``average
out'' over some bounded subset of the phase space, instead of blowing
up to infinity (as a simplest example, one could think of a periodic
dependence on $s$). In Appendix \ref{sec:explicit_time} we set up a
suitable mathematical framework, in which the above arguments are
arranged in a more formal fashion, through a combination of the
appropriate skew product flows (see, for example,
\cite{CarLukRea,CraDebFla} and references therein). The resulting
averaged system is given by
\begin{subequations}
\label{eq:explicit_time_averaged}
\begin{equation}
\dif\BS x=\langle\BS F\rangle(\BS x,t)\dif t+
\langle\BS\kappa\rangle(\BS x,t)\dif\BS V_t,
\end{equation}
\begin{equation}
\langle\BS F\rangle(\BS x,t)=\lim_{r\to\infty}\frac 1r\int_0^r
\BS F(\BS x,\BS z(s),t,s)\dif s
\end{equation}
\begin{equation}
\begin{split}
\langle\BS\kappa\rangle(\BS x,t)\langle\BS\kappa\rangle(\BS x,t)^T
=\lim_{r\to\infty}\frac 1r\int_0^r\BS\kappa(\BS x,t,s)\BS\kappa(\BS
x,t,s)^T\dif s.
\end{split}
\end{equation}
\begin{equation}
\dif\BS{TX}_{\BS x}^{t_0,t}=\left[\left(\left\langle\parderiv{\BS F}
{\BS x}\right\rangle(\BS x,t)+\BS H(\BS x,t)\right)\dif t+
\parderiv{\langle\BS\kappa\rangle(\BS x,t)}{\BS x}\dif\BS V_t\right]
\BS{TX}_{\BS x}^{t_0,t},
\end{equation}
\begin{equation}
\left\langle\parderiv{\BS F}{\BS x}\right\rangle(\BS x,t)=
\lim_{r\to\infty}\frac 1r\int_0^r\parderiv{\BS F}{\BS x}(\BS x,\BS z(s)
,t,s)\dif s,
\end{equation}
\begin{equation}
\label{eq:H_corr_nonauto}
\begin{split}
\BS H(\BS x,t)=\int_0^\infty\bigg(\lim_{r\to\infty}\frac 1r
\int_0^r\parderiv{\BS F}{\BS y}(\BS x,\BS z(\tau+s),t,\tau+s)
\times\\\times \BS T_{\BS x,t,\BS z(\tau)}^{\tau,s}\parderiv{\BS
  G}{\BS x}(\BS x,\BS z(\tau),t,\tau)\dif\tau\bigg)\dif s,
\end{split}
\end{equation}
\end{subequations}
where $\BS T_{\BS x,t,\BS z(\tau)}^{\tau,s}$ is the tangent map of
\eqref{eq:dyn_sys_nonauto_fast_limiting_z} at $\BS z(\tau)$ to
$\tau+s$.

Now, following Section \ref{sec:coupling}, we write
\eqref{eq:dyn_sys_nonauto} in the linear energy-preserving coupling
form:
\begin{subequations}
\label{eq:dyn_sys_nonauto_en}
\begin{equation}
\label{eq:dyn_sys_nonauto_slow_en}
\dif\BS x=\left[\BS f(\BS x,t,t/\varepsilon)+\BS S_x^{-1/2}\BS L(t)
\BS S_y^{1/2}\BS y\right]\dif t+\BS\kappa(\BS x,t,t/\varepsilon)
\dif\BS V_t,
\end{equation}
\begin{equation}
\label{eq:dyn_sys_nonauto_fast_en}
\dif\BS y=\frac 1\varepsilon\left[\BS g(\BS y,t,t/\varepsilon)-
\BS S_y^{-1/2}\BS L^T(t)\BS S_x^{1/2}\BS x\right]\dif t+
\frac 1{\sqrt\varepsilon}\BS\sigma(\BS y,t,t/\varepsilon)\dif\BS W_t.
\end{equation}
\end{subequations}
Here, observe that the linear coupling matrix $\BS L$ may depend on
the slow time $t$, $\BS L=\BS L(t)$. With
\eqref{eq:dyn_sys_nonauto_en} it is not hard to see that that averaged
dynamics for the slow tangent map are given by
\begin{subequations}
\label{eq:dyn_sys_nonauto_slow_tangent_en}
\begin{equation}
\begin{split}
\dif\BS{TX}_{\BS x}^{t_0,t}=\bigg[\bigg(\bigg\langle\parderiv{\BS f}
{\BS x}\bigg\rangle(\BS x,t)-\BS S_x^{-1/2}\BS L(t)\BS S_y^{1/2}
\BS{\mathcal C}(\BS x,t)\BS S_y^{-1/2}\BS L^T(t)\BS S_x^{1/2}
\bigg)\dif t+\\+\parderiv{\langle\BS\kappa\rangle(\BS x,t)}{\BS x}
\dif\BS V_t\bigg]\BS{TX}_{\BS x}^{t_0,t},
\end{split}
\end{equation}
\begin{equation}
\left\langle\parderiv{\BS f}{\BS x}\right\rangle(\BS x,t)=
\lim_{r\to\infty}\frac 1r\int_0^r\parderiv{\BS f}{\BS x}(\BS x,t,s)
\dif s,
\end{equation}
\begin{equation}
\begin{split}
\BS{\mathcal C}(\BS x,t)=\int_0^\infty\left(\lim_{r\to\infty}\frac
1r\int_0^r \BS T_{\BS x,t,\BS z(\tau)}^{\tau,s}\dif\tau\right)\dif s.
\end{split}
\end{equation}
\end{subequations}
Here we can see that the situation is completely analogous to what was
found in Sections \ref{sec:coupling}--\ref{sec:suppression}, that is,
the positive-definiteness of $\BS{\mathcal C}(\BS x,t)$ in the
canonical energy coordinates has a damping effect on the slow averaged
tangent dynamics in \eqref{eq:dyn_sys_nonauto_slow_tangent_en},
thereby reducing chaos at the slow variables.

\section{Nonlinear energy-preserving coupling with dependence on the
slow and fast variables}
\label{sec:slow_fast_dep}

Recall that in Sections \ref{sec:coupling} and \ref{sec:extension} we
assumed that the matrix $\BS L$ in the energy preserving coupling is
either a constant matrix (Section \ref{sec:coupling}), or is a
function of only the slow time $t$ (Section
\ref{sec:extension}). However, for the purpose of the energy
preservation, $\BS L$ can be an arbitrary $N_x\times N_y$ matrix
function $\BS L(\BS x,\BS y,t,t/\varepsilon)$. Below we briefly
outline the situation when $\BS L$ is a function of $\BS x$, $\BS y$,
$t$ and $t/\varepsilon$.

In such a case, we can extract the slow part $\BS L(t)$ of the
energy-preserving coupling $\BS L(\BS x,\BS y,t,t/\varepsilon)$ by
time-averaging over the fast variables on the trajectory $\BS x(t)$:
\begin{subequations}
\label{eq:L_decomp}
\begin{equation}
\BS L(t)=\lim_{r\to\infty}\frac 1r\int_0^r\BS L(\BS x(t),\BS
z(s),t,s)\dif s,
\end{equation}
\begin{equation}
\BS L(\BS x,\BS y,t,s)=\BS L(t)+\BS L^\prime(\BS x,\BS y,t,s),
\end{equation}
\end{subequations}
where $\BS x(t)$ is the solution of \eqref{eq:dyn_sys_nonauto_en}, and
$\BS z(s)$ is the solution of
\eqref{eq:dyn_sys_nonauto_fast_limiting_z} for given $\BS x$ and
$t$. Thus, the term $\BS L^\prime(\BS x,\BS y,t,s)$ in
\eqref{eq:L_decomp} represents rapid oscillations with zero mean state
around the slowly varying $\BS L(t)$:
\begin{equation}
\lim_{r\to\infty}\frac 1r\int_0^r\BS L^\prime(\BS x(t),\BS
z(s),t,s)\dif s =\BS 0.
\end{equation}
With \eqref{eq:L_decomp}, the equations in \eqref{eq:dyn_sys_nonauto_en}
become
\begin{subequations}
\begin{equation}
\dif\BS x=\left[\BS f(\BS x,\BS y,t,t/\varepsilon)+\BS
S_x^{-1/2}\BS L(t)\BS S_y^{1/2}\BS y\right]\dif t+
\BS\kappa(\BS x,t,t/\varepsilon)\dif\BS V_t,
\end{equation}
\begin{equation}
\dif\BS y=\frac 1\varepsilon\left[\BS g(\BS x,\BS y,t,t/\varepsilon)-
\BS S_y^{-1/2}\BS L^T(t)\BS S_x^{1/2}\BS x\right]\dif t+\frac
1{\sqrt\varepsilon}\BS\sigma(\BS y,t,t/\varepsilon)\dif\BS W_t,
\end{equation}
\end{subequations}
where we denote
\begin{subequations}
\begin{equation}
\BS f(\BS x,\BS y,t,s)=\BS f(\BS x,t,s)+\BS S_x^{-1/2}\BS L^\prime
(\BS x,\BS y,t,s)\BS S_y^{1/2}\BS y,
\end{equation}
\begin{equation}
\BS g(\BS x,\BS y,t,s)=\BS g(\BS y,t,s)-\BS S_y^{-1/2}\BS L^{\prime
  T}(\BS x,\BS y,t,s) \BS S_x^{1/2}\BS x,
\end{equation}
\end{subequations}
With the above changes, the right-hand side of
\eqref{eq:dyn_sys_nonauto_slow_tangent} acquires a new term:
\begin{subequations}
\begin{equation}
\begin{split}
\dif\BS{TX}^{t_0,t}_{\BS x}=\bigg[\bigg(\left\langle\parderiv{\BS f}{\BS x}
\right\rangle(\BS x,t)-\BS S_x^{-1/2}\BS L(t)\BS S_y^{1/2}
\BS{\mathcal C}(\BS x,t)\BS S_y^{-1/2}\BS L^T(t)\BS S_x^{1/2}+\\+
\int_0^\infty\BS H'(\BS x,t,s)\bigg)\dif s+
\parderiv{\langle\BS\kappa\rangle(\BS x,t)}{\BS x}
\dif\BS V_t\bigg]\BS{TX}^{t_0,t}_{\BS x},
\end{split}
\end{equation}
\begin{equation}
\begin{split}
\BS H^\prime(\BS x,t,s)=\BS S_x^{-1/2}\BS L(t)\BS S_y^{1/2}
\left[\lim_{r\to\infty}\frac 1r\int_0^r\BS T_{\BS x,t,\BS z(\tau)}^{\tau,s}
  \parderiv{\BS g}{\BS x}(\BS x,\BS
  z(\tau),t,\tau)\dif\tau\right]-\\-\left[\lim_{r\to\infty}\frac
  1r\int_0^r\parderiv{\BS f}{\BS z} (\BS x,\BS z(\tau+s),t,\tau+s)\BS
  T_{\BS x,t,\BS z(\tau)}^{\tau,s}\dif\tau\right] \BS S_y^{-1/2}\BS L^T(t)\BS
S_x^{1/2}+\\+\lim_{r\to\infty} \frac1r\int_0^r\parderiv{\BS f}{\BS
  z}(\BS x,\BS z(\tau+s),t,\tau+s)\BS T_{\BS x,t,\BS
  z(\tau)}^{\tau,s}\parderiv{\BS g}{\BS x}(\BS x,\BS z(\tau),t,\tau)\dif\tau.
\end{split}
\end{equation}
\end{subequations}
It is not known at present how the new term $\BS H^\prime(\BS x,t,s)$
is going to affect the dynamics, however, the chaos-suppressing term
$\BS{\mathcal C}(\BS x,t)$ is nonetheless retained by the dynamics.
One can presume that if the effect of the new term is small, the chaos
suppression through $\BS{\mathcal C}(\BS x,t)$ would still be
observed. At this point, it is hardly possible to say anything about
the effect of $\BS H^\prime(\BS x,t,s)$ without making additional
assumptions about its structure. Clearly, the answer to this question
depends on a particular type of energy-preserving coupling used in the
model. While a general coupling of the form in \eqref{eq:L_decomp}
might not necessarily promote the suppression of chaos at the slow
variables, there could be additional restrictions on $\BS L^\prime(\BS
x,\BS y,t,s)$ imposed by the physics of a particular model under
consideration. In the future work, the author intends to derive
suitable classes of a more general nonlinear energy-preserving
coupling which nonetheless promotes the suppression of chaos at the
slow variables.

\section{Invariance with respect to time-scaling and independence
of $\varepsilon$}
\label{sec:invariance}

Throughout the work, it was presumed that the dependence of the model
on the time-scale separation parameter $\varepsilon$ is explicit.
However, often in real-world situations the parameter $\varepsilon$ is
not present explicitly, but it is known from the observations that
there is a subset of the ``slow'' variables in the system, while the
rest of the variables are ``fast''. It turns out that the framework,
developed here, also does not explicitly depend on
$\varepsilon$. Observe that the averages in
\eqref{eq:dyn_sys_nonauto_slow_tangent} are invariant under the
time-rescaling of the limiting dynamics in
\eqref{eq:dyn_sys_nonauto_fast_limiting_z} as
\begin{equation}
\label{eq:dyn_sys_fast_limiting_z_alpha}
\dif\BS z=\alpha\BS G(\BS x,\BS z,t,s)\dif s+ \sqrt\alpha\,
\BS\sigma(\BS z,t,s)\dif\BS W_s,
\end{equation}
where $\alpha>0$ is an arbitrary parameter. Setting
$\alpha=\varepsilon^{-1}$ yields
\begin{equation}
\label{eq:dyn_sys_fast_limiting_z_eps}
\dif\BS z=\frac 1\varepsilon\BS G(\BS x,\BS z,t,s)\dif s+ \frac
1{\sqrt\varepsilon}\BS\sigma(\BS z,t,s)\dif\BS W_s,
\end{equation}
which is obtained directly from the fast equation in
\eqref{eq:dyn_sys_nonauto_fast} by setting $\BS x$ and $t$ as constant
parameters. Due to this property, the explicit knowledge of the
time-scale separation parameter $\varepsilon$ is not necessary; for
all practical purposes, one can assume that the original slow-fast
system in \eqref{eq:dyn_sys_nonauto} is supplied without the
time-scale separation parameter $\varepsilon$, in the form
\begin{subequations}
\label{eq:dyn_sys_no_eps}
\begin{equation}
\label{eq:dyn_sys_no_eps_slow}
\dif\BS x=\BS F(\BS x,\BS y,t,t)\dif t+\BS\kappa(\BS x,t,t)\dif\BS V_t,
\end{equation}\
\begin{equation}
\label{eq:dyn_sys_no_eps_fast}
\dif\BS y=\BS G(\BS x,\BS y,t,t)\dif t+\BS\sigma(\BS y,t,t)\dif\BS W_t,
\end{equation}
\end{subequations}
for which it is {\em known}, that $\BS x$ is the set of slow variables
and $\BS y$ is the set of fast variables, and the dependence on the
first $t$-parameter in each function is slow, while the dependence on
the second $t$-parameter is fast. The limiting dynamics in
\eqref{eq:dyn_sys_nonauto_fast_limiting_z} are thus obtained by
setting $\BS x$ and first $t$-parameter in
\eqref{eq:dyn_sys_no_eps_fast} to constants, and the averaged slow
tangent dynamics in \eqref{eq:dyn_sys_nonauto_slow_tangent} remain in
effect with no corrections.

\section{Summary}
\label{sec:summary}

In this work we study the effect of the fast rapidly mixing dynamics
on the chaos at slow variables in a two-scale system with a
deterministic linear energy-preserving coupling. A suitable theory is
developed by applying the averaging formalism to the tangent dynamics
of the system, and it is found that the linear energy-preserving
coupling creates a systematic damping effect on the chaos at the slow
scales when the fast dynamics is rapidly mixing (the precise criterion
is the positive-definiteness of the infinite-time linear response to a
small constant external forcing at the fast scales), which can also be
interpreted as the linear stability property of the long-term mean
state of the fast limiting system under the external perturbations. It
is suggested that this property must strongly depend on the ratio of
the autocorrelation time and the advection time scale at the fast
variables, however, a systematic study of this proposition is left for
the future work. The effect of the chaos suppression at the slow
variables through the rapidly mixing fast dynamics is systematically
demonstrated for the two-scale Lorenz model, which is also
appropriately rescaled so that the adjustments for the mixing regime
at the fast variables do not affect the mean state and variance of
both the slow and fast variables. In particular, it is shown through
the numerical experiment that the uncoupled slow dynamics, where the
fast variables are replaced with their averages (which are zeros in
the rescaled Lorenz model), may remain chaotic, while the full coupled
system loses chaos and becomes completely predictable at the slow
scales as the dynamics at the fast scales become more turbulent. The
theory is also systematically generalized onto the more general
setting with nonautonomous dynamics and stochastic noise at both the
slow and fast variables.

\subsection*{Future work}

Clearly, the present work, while uncovering a dynamical mechanism for
the chaos suppression at the slow variables in a multiscale system,
does not yet provide any systematic quantification of the chaos
suppression effect, which will have to be addressed in the future
work. Below we sketch a few key directions, immediately emerging from
the results of the present study.
\begin{itemize}
\item {\bf Is there a simple direct quantitative connection between
  the chaos at the slow variables and the infinite-time linear
  response at the fast variables?} Above, while the dynamical
  mechanism of the chaos suppression is discussed, no quantitative
  connections between the Lyapunov characteristic exponents of the
  slow averaged dynamics and the eigenvalues of the symmetric part of
  the infinite-time linear response at the fast variables are
  investigated. Are there any simple direct estimates of the chaos
  reduction at the slow variables through the autocorrelation times of
  the fast variables?
\item {\bf Is it possible to estimate the chaos suppression effect
  from the autocorrelation time/advection time scale ratio?} Above we
  briefly formulated the hypothesis that the suppression of chaos at
  the slow variables could be amplified by a large ratio of the
  advection time scale to the autocorrelation time at the fast
  variables. If this hypothesis is correct, would it be possible to
  make an estimate of the chaos suppression effect by knowing the
  ratio of the advection time scale to the autocorrelation time,
  perhaps through the theory of turbulence?
\item {\bf What if the energy-preserving coupling is nonlinear?} We
  briefly sketched that if the energy-preserving coupling is
  nonlinear, then an additional term appears in the dynamics of the
  slow averaged tangent map, and it is not immediately clear whether
  it would create a damping or amplifying slow chaos effect. Can one
  identify suitable physically relevant classes of nonlinear coupling
  which suppresses chaos at the slow variables?
\item {\bf What happens in more sophisticated models?}  Are the fully
  coupled models less chaotic at the slow variables than the uncoupled
  low-frequency variability models? If so, could the global climate be
  generally less chaotic than we presently tend to think? Could the
  predictive skill be improved by appropriately adjusting the uncoupled
  models to match the slow chaos properties of the coupled dynamics?
\end{itemize}
The author intends to address these problems in the future work,
possibly in collaboration with atmospheric scientists.

\medskip

{\bf Acknowledgment.}
\medskip

The author thanks Ibrahim Fatkullin for fruitful discussions. The
author is supported by the National Science Foundation CAREER grant
DMS-0845760, and the Office of Naval Research grants N00014-09-0083
and 25-74200-F6607.

\appendix

\section{Infinite time linear response to changes in $\BS x$}
\label{sec:inf_time_lin_resp}

Let $\mu_{\BS x}$ be an invariant measure for the dynamics given by
\begin{equation}
\label{eq:fast_lim}
\deriv{\BS z}t=\BS G(\BS x,\BS z)
\end{equation}
with flow $\phi_{\BS x}^t$ and ergodic attractor $\mathcal A_{\BS x}$,
where $\BS x$ is a fixed parameter. Let $B(\BS z)$ be a differentiable
function of $\BS z$, and let $\langle B\rangle_{\BS x}$ denote the
average
\begin{equation}
\langle B\rangle_{\BS x}=\int_{\mathcal A_{\BS x}}B(\BS z)\dif\mu_{\BS
  x}(\BS z).
\end{equation}
Our goal here is to compute the derivative of $\langle B\rangle_{\BS
  x}$ with respect to $\BS x$ (if it exists):
\begin{equation}
\label{eq:B_deriv}
\deriv{\langle B\rangle_{\BS x}}{\BS x}=\int_{\mathcal A_{\BS x}}B(\BS
z)\dif\mu_{\BS x}^\prime(\BS z).
\end{equation}
Here we are going to use the linear response theory from
\cite{Abr5,Abr6,Abr7,AbrMaj4,AbrMaj5,AbrMaj6,MajAbrGro,Rue2,Ris}. Let
us assume that, at the initial time $t=0$, there is a statistical
ensemble of points $\BS z$, distributed according to $\mu_{\BS x}$,
and, at the same time, there is a change $\BS x+\delta\BS x$ in
\eqref{eq:fast_lim}. Then, one can treat the difference between the
averages along the perturbed and unperturbed trajectories of
\eqref{eq:fast_lim}, emerging from the initial statistical ensemble,
as time goes to infinity, as the derivative in \eqref{eq:B_deriv}:
\begin{equation}
\label{eq:B_deriv_2}
\begin{split}
\int_{\mathcal A_{\BS x}}&B(\BS z)\dif\mu_{\BS x}^\prime(\BS
z)\delta\BS x+O(\|\delta\BS x\|^2)=\\=&\lim_{t\to\infty}\int_{\mathcal A_{\BS x}}
\left[B(\phi_{\BS x+\delta\BS x}^t\BS z)-B(\phi_{\BS x}^t\BS z)\right]
\dif\mu_{\BS x}(\BS z).
\end{split}
\end{equation}
For the above relation to hold, the structural stability of $\mu_{\BS
  x}$ under the changes in $\BS x$ is required. It is known that
uniformly hyperbolic diffeomorphisms on the whole phase space of
\eqref{eq:fast_lim} (Anosov) or just its nonwandering set (Axiom A)
are structurally stable \cite{EckRue,Rue,Rue2,Rue3,You}. The integral
in the right-hand side above is given by
\begin{equation}
\label{eq:app1_temp}
\begin{split}
\int_{\mathcal A_{\BS x}}&\left[B(\phi_{\BS x+\delta\BS x}^t\BS z)-
B(\phi_{\BS x}^t\BS z)\right]\dif\mu_{\BS x}(\BS z)=\\=&
\int_{\mathcal A_{\BS x}}\left[\nabla B(\phi_{\BS x}^t\BS z)\parderiv{}{\BS x}
\phi_{\BS x}^t\BS z\right]\dif\mu_{\BS x}(\BS z)\delta\BS x+
O(\|\delta\BS x\|^2).
\end{split}
\end{equation}
To compute the derivative of the flow $\phi_{\BS x}^t\BS z$ in $\BS
x$, we use the fact that it is a solution of \eqref{eq:fast_lim}, and,
therefore, satisfies
\begin{equation}
\parderiv{}t\phi_{\BS x}^t\BS z=\BS G(\BS x,\phi_{\BS x}^t\BS z),\qquad
\parderiv{}t\phi_{\BS x+\delta\BS x}^t\BS z=\BS G(\BS x+\delta\BS x,
\phi_{\BS x+\delta\BS x}^t\BS z).
\end{equation}
Subtracting one identity from another and Taylor-expanding in
$\delta\BS x$, we obtain
\begin{equation}
\parderiv{}t\left(\parderiv{}{\BS x}\phi_{\BS x}^t\BS z\right)=
\parderiv{\BS G}{\BS z}(\BS x,\phi_{\BS x}^t\BS
z)\left(\parderiv{}{\BS x} \phi_{\BS x}^t\BS z\right)+\parderiv{\BS
  G}{\BS x}(\BS x,\phi_{\BS x}^t\BS z)+O(\|\delta\BS x\|).
\end{equation}
By sending $\delta\BS x$ to zero, the formal solution to the above
equation with zero initial condition is given by the Duhamel's
principle:
\begin{equation}
\parderiv{}{\BS x}\phi_{\BS x}^t\BS z=\int_0^t\BS T_{\BS x,\BS
  z}^{t-s} \parderiv{\BS G}{\BS x}(\BS x,\phi_{\BS x}^s\BS
z)\dif s,
\end{equation}
where $\BS T_{\BS x,\BS z}^t$ is the tangent map of $\phi_{\BS x}^t\BS
z$
\begin{equation}
\BS T_{\BS x,\BS z}^t=\parderiv{}{\BS z}\phi_{\BS x}^t\BS z,
\end{equation}
satisfying
\begin{equation}
\parderiv{}t\BS T_{\BS x,\BS z}^t=\parderiv{\BS G}{\BS z}(\BS x,
\phi_{\BS x}^\tau\BS z)\BS T_{\BS x,\BS z}^t.
\end{equation}
Substituting the above equation into \eqref{eq:app1_temp}, we obtain
\begin{equation}
\begin{split}
\int_{\mathcal A_{\BS x}}&\left[B(\phi_{\BS x+\delta\BS x}^t\BS z)-
  B(\phi_{\BS x}^t\BS z)\right]\dif\mu_{\BS x}(\BS z)=\\&=
\int_0^t\int_{\mathcal A_{\BS x}}\left[\nabla B(\phi_{\BS x}^t\BS z)
  \BS T_{\BS x,\BS z}^{t-s} \parderiv{\BS G}{\BS x}(\BS x,\phi_{\BS
    x}^s\BS z)\right]\dif\mu_{\BS x}(\BS z)\dif s\delta\BS x+
O(\|\delta\BS x\|^2).
\end{split}
\end{equation}
By using the fact that $\mu_{\BS x}$ is the invariant measure for
$\phi_{\BS x}^t$, we obtain
\begin{equation}
\begin{split}
\int_{\mathcal A_{\BS x}}&\left[B(\phi_{\BS x+\delta\BS x}^t\BS z)-
  B(\phi_{\BS x}^t\BS z)\right]\dif\mu_{\BS x}(\BS z)=\\&=
\int_0^t\int_{\mathcal A_{\BS x}}\left[\nabla B(\phi_{\BS
    x}^{t-s}\BS z) \BS T_{\BS x,\BS z}^{t-s} \parderiv{\BS
    G}{\BS x}(\BS x,\BS z)\right] \dif\mu_{\BS x}(\BS
z)\dif s\delta\BS x+ O(\|\delta\BS x\|^2)=\\&=
\int_0^t\int_{\mathcal A_{\BS x}}\left[\nabla B(\phi_{\BS
    x}^s\BS z) \BS T_{\BS x,\BS z}^s\parderiv{\BS
    G}{\BS x}(\BS x,\BS z)\right] \dif\mu_{\BS x}(\BS
z)\dif s\delta\BS x+ O(\|\delta\BS x\|^2).
\end{split}
\end{equation}
Comparing the above relation with \eqref{eq:B_deriv_2}, we obtain, as
$\delta\BS x\to 0$,
\begin{equation}
\int_{\mathcal A_{\BS x}}B(\BS z)\dif\mu_{\BS x}^\prime(\BS z)=
\int_0^\infty\int_{\mathcal A_{\BS x}}\left[\nabla B(\phi_{\BS x}^s
\BS z) \BS T_{\BS x,\BS z}^s\parderiv{\BS G}{\BS x}(\BS x,\BS z)
\right] \dif\mu_{\BS x}(\BS z)\dif s.
\end{equation}
By using ergodicity of $\mu_{\BS x}$, we can replace the measure
averages above with time averages over a long-term trajectory $\BS
z(t)$ of \eqref{eq:fast_lim}:
\begin{equation}
\int_{\mathcal A_{\BS x}}B(\BS z)\dif\mu_{\BS x}^\prime(\BS z)=
\int_0^\infty\left(\lim_{r\to\infty}\frac 1r\int_0^r\left[\nabla
B(\BS z(\tau+s))\BS T_{\BS x,\BS z(\tau)}^s\parderiv{\BS G}{\BS x}
(\BS x,\BS z(\tau))\right]\dif\tau\right)\dif s.
\end{equation}
In terms of the linear response theory
\cite{Abr5,Abr6,Abr7,AbrMaj4,AbrMaj5,AbrMaj6,MajAbrGro,Rue2,Ris}, the
above expression is the infinite time linear response of $\langle
B\rangle$ to the changes in $\BS x$.

\section{Infinite time linear response of the mean state to
constant additive perturbation}
\label{sec:inf_time_add}

Let $\mu$ be an invariant measure for the dynamics given by
\begin{equation}
\label{eq:fast_lim_add}
\deriv{\BS z}t=\BS G(\BS z)
\end{equation}
with flow $\phi^t$ and ergodic attractor $\mathcal A$. Let $\langle\BS
z\rangle$ be the mean state of \eqref{eq:fast_lim_add}:
\begin{equation}
\langle\BS z\rangle=\int_{\mathcal A}\BS z\dif\mu(\BS z).
\end{equation}
Now, assume that \eqref{eq:fast_lim_add} is perturbed by a constant
additive perturbation $\delta\BS f$,
\begin{equation}
\label{eq:fast_lim_add_pert}
\deriv{\BS z}t=\BS G(\BS z)+\delta\BS f,
\end{equation}
with perturbed flow given by $\phi^{*t}$. Our goal here is to compute
the infinite time linear response of $\langle\BS z\rangle$. Following
Appendix \ref{sec:inf_time_lin_resp}, we write
\begin{equation}
\parderiv{}t\phi^t\BS z=\BS G(\phi^t\BS z),\qquad
\parderiv{}t\phi^{*t}\BS z=\BS G(\phi^{*t}\BS z)+\delta\BS f.
\end{equation}
Let the difference between the perturbed and unperturbed flows be
given as $\delta\phi^t\BS z$. Now, taking the difference between the
two relations above, we obtain
\begin{equation}
\begin{split}
\parderiv{}t(\delta\phi^t\BS z)=\BS G(\phi^t\BS z+\delta\phi^t\BS z)-
\BS G(\phi^t\BS z)+\delta\BS f=\\=\nabla\BS G(\phi^t\BS z)
\delta\phi^t\BS z+\delta\BS f+O(\|\delta\phi^t\BS z\|^2).
\end{split}
\end{equation}
Neglecting the higher-order term, we obtain, through the Duhamel's
principle,
\begin{equation}
\delta\phi^t\BS z=\int_0^t\BS T_{\BS z}^s\dif s\,\delta\BS f,
\end{equation}
where $\BS T_{\BS z}^t$ is the tangent map of \eqref{eq:fast_lim_add},
given by
\begin{equation}
\parderiv{}t\BS T_{\BS z}^t=\nabla\BS G(\phi^t\BS z)\BS T_{\BS z}^t.
\end{equation}
Then, by integrating over $\mu$ and taking $t\to\infty$, for the
infinite time response of $\langle\BS z\rangle$ we obtain
\begin{equation}
\delta\langle\BS z\rangle=\int_0^\infty\int_{\mathcal A}\BS T_{\BS z}^s
\dif\mu(\BS z)\dif s\,\delta\BS f.
\end{equation}
Assuming ergodicity of $\mu$ and replacing measure average with time
average over a long-term trajectory $\BS z(t)$ of
\eqref{eq:fast_lim_add}, we obtain
\begin{equation}
\delta\langle\BS z\rangle=\int_0^\infty\left(\lim_{r\to\infty}\frac 1r
\int_0^r\BS T_{\BS z(\tau)}^s\dif\tau\right)\dif s\,\delta\BS f.
\end{equation}

\section{Averaged dynamics with explicit time dependence}
\label{sec:explicit_time}

Here, we assume that the fast limiting system in
\eqref{eq:dyn_sys_nonauto_fast_limiting_z}, which includes both the
$s$-dependent Wiener process $\BS W_s$ and the explicit $s$-dependence
of $\BS G$ and $\BS\sigma$, can be expressed as a collection of the
following dynamical systems:
\begin{itemize}
\item A flow $\theta^s:\mathcal Q\to\mathcal Q$ on a measurable space
  $\mathcal Q$ with a distribution measure $\rho$, where $s$ is the
  time parameter. The elements $q\in\mathcal Q$ play the role of a
  bounded pseudo-``time'' in the deterministic time-dependence of $\BS
  F$, $\BS G$, $\BS\kappa$ and $\BS\sigma$. $\theta^s$ has the
  following properties:
  \begin{itemize}
    \item The semigroup property -- $\theta^0=I$,
      $\theta^s\theta^r=\theta^{s+r}$;
    \item $\theta$ is measure-preserving: $\rho(Q)=\rho(\theta^sQ)$
      for any subset $Q\subset\mathcal Q$.
  \end{itemize}
\item A base flow $\vartheta^s:\Omega\to\Omega$, given over a
  probability space $\{\Omega,\mathcal F,\mathbb P\}$ (this
  probability space documents all realizations of the Wiener process
  $\BS W_s(\omega)$). The idea of this dynamical system is to
  ``fast-forward'' or ``rewind'' $\BS W_s(\omega)$ as needed by
  specifying the time parameter $s$.  $\vartheta^s$ has the following
  properties:
  \begin{itemize}
    \item The semigroup property -- $\vartheta^0=I$,
      $\vartheta^s\vartheta^r=\vartheta^{s+r}$;
    \item If $A\in\mathcal F$, then $\vartheta^sA$ is also in
      $\mathcal F$ for any $s$;
    \item Measure preservation: $\mathbb P(A)=\mathbb P(\vartheta^sA)$
      for any $s$;
    \item The Wiener process $\BS W_s(\omega)$ satisfies $\BS
      W_r(\vartheta^s\omega)=\BS W_{s+r}(\omega)-\BS W_s(\omega)$.
  \end{itemize}
\item A $\theta,\vartheta$-cocycle $\phi_{\BS
  x,t}^{q,\omega,s}:\mathbb R^{N_y}\to\mathbb R^{N_y}$ with properties
\begin{equation}
\label{eq:cocycle}
\phi_{\BS x,t}^{q,\omega,0}=I,\qquad \phi_{\BS
  x,t}^{q,\omega,s+r}=\phi_{\BS x,t}^{\theta^r
  q,\vartheta^r\omega,s}\;\;\phi_{\BS x,t}^{q,\omega,r},\quad
\forall(q,\omega)\in(\mathcal Q,\Omega),
\end{equation}
which represents a solution of the non-autonomous fast limiting system
in \eqref{eq:dyn_sys_nonauto_fast_limiting_z}.
\end{itemize}
Now, instead of considering a still compact global attractor $\mathcal
A_{\BS x}$ for the autonomous dynamics like in Section
\ref{sec:twoscale} (which parametrically depends on the slow variables
$\BS x$), here we consider a $q,\omega$-dependent family of compact
global pullback-attracting sets $\mathcal A_{\BS
  x,t}(q,\omega)\subset\mathbb R^{N_y}$, which is
$\theta,\vartheta$-invariant, that is,
\begin{equation}
\phi_{\BS x,t}^{q,\omega,s}\mathcal A_{\BS x,t}(q,\omega)=\mathcal
A_{\BS x,t} (\theta^sq,\vartheta^s\omega),\qquad
\lim_{s\to\infty}dist(\mathcal A_{\BS x,t}(q,\omega),
\phi_{\BS x,t}^{\theta^{-s}q,\vartheta^{-s}\omega,s}\mathcal B)=0,
\end{equation}
for all $q\in\mathcal Q$, $\mathbb P$-almost all $\omega\in\Omega$,
and all compact subsets $\mathcal B$ of $\mathbb R^{N_y}$, where
$dist$ is the Hausdorff distance in $\mathbb R^{N_y}$. Then, the
family of sets $\mathcal A_{\BS x,t}(q,\omega)$ is called the global
pullback attractor of
\eqref{eq:dyn_sys_nonauto_fast_limiting_z}. Similarly, instead of a
single invariant measure $\mu_{\BS x}$ one can consider a family of
$\theta,\vartheta$-invariant measures $\mu_{\BS x,t}^{q,\omega}$ on
$\mathcal A_{\BS x,t}(q,\omega)$, with the property
\begin{equation}
\label{eq:mu_inv}
\mu_{\BS x,t}^{q,\omega}(\phi_{\BS x,t}^{q,\omega,s}A)=\mu_{\BS
  x,t}^{\theta^sq,\vartheta^s\omega}(A), \qquad\forall
A\subset\mathcal A_{\BS x,t}(q,\omega).
\end{equation}
At this point, let
\begin{equation}
\begin{split}
\BS{\tilde F}(\BS x,\BS z,t,\theta^sq_0)=\BS F(\BS x,\BS z,t,s),&\qquad
\BS{\tilde G}(\BS x,\BS z,t,\theta^sq_0)=\BS G(\BS x,\BS z,t,s),\\
\BS{\tilde\kappa}(\BS x,t,\theta^sq_0)=\BS\kappa(\BS x,t,s),&\qquad
\BS{\tilde\sigma}(\BS z,t,\theta^sq_0)=\BS\sigma(\BS z,t,s),
\end{split}
\end{equation}
for all $\BS x$, $t$, $\BS z$, $s$ and some $q_0\in\mathcal Q$, which
corresponds to the zero time $s=0$. Now, the slow limiting averaged
system for \eqref{eq:dyn_sys_nonauto} is given by
\begin{subequations}
\label{eq:dyn_sys_nonauto_slow_limiting_x}
\begin{equation}
\dif\BS x=\langle\BS F\rangle(\BS x,t)\dif t+
\langle\BS\kappa\rangle(\BS x,t)\dif\BS V_t,
\end{equation}
\begin{equation}
\langle\BS F\rangle(\BS x,t)=\int_{\mathcal Q}\int_{\Omega}
\int_{\mathcal A_{\BS x,t}(q,\omega)}\BS{\tilde F}(\BS x,\BS z,t,q)
\dif\mu_{\BS x,t}^{q,\omega}(\BS z)\dif\mathbb P(\omega)\dif\rho(q),
\end{equation}
\begin{equation}
\begin{split}
\langle\BS\kappa\rangle(\BS x,t)\langle\BS\kappa\rangle(\BS x,t)^T=
\int_{\mathcal Q}\BS{\tilde\kappa}(\BS x,t,q)\BS{\tilde\kappa}(\BS
x,t,q)^T\dif\rho(q)=\\=\lim_{r\to\infty}\frac 1r\int_0^r\BS\kappa(\BS
x,t,s)\BS\kappa(\BS x,t,s)^T\dif s.
\end{split}
\end{equation}
\end{subequations}
Here we assume that the cocycle in \eqref{eq:cocycle}, which is a
stochastic flow in $s$, is differentiable with respect to its initial
condition \cite{Kun}, and, therefore, the corresponding equations for
the averaged slow tangent map in \eqref{eq:dyn_sys_slow_tangent} are
given by
\begin{subequations}
\label{eq:dyn_sys_nonauto_slow_tangent}
\begin{equation}
\dif\BS{TX}_{\BS x}^{t_0,t}=\left[\left(\left\langle\parderiv{\BS F}
{\BS x}\right\rangle(\BS x,t)+\BS H(\BS x,t)\right)\dif t+
\parderiv{\langle\BS\kappa\rangle(\BS x,t)}{\BS x}\dif\BS V_t\right]
\BS{TX}_{\BS x}^{t_0,t},
\end{equation}
\begin{equation}
\left\langle\parderiv{\BS F}{\BS x}\right\rangle(\BS x,t)=
\int_{\mathcal Q}\int_\Omega\int_{\mathcal A_{\BS x,t}(q,\omega)}
\parderiv{\BS{\tilde F}} {\BS x}(\BS x,\BS z,t,q)\dif\mu_{\BS x,t}
^{q,\omega}(\BS z)\dif\mathbb P(\omega)\dif\rho(q),
\end{equation}
\begin{equation}
\BS H(\BS x,t)=\int_{\mathcal Q}\int_\Omega\int_{\mathcal A_{\BS x,t}(q,\omega)}
\BS{\tilde F}(\BS x,\BS z,t,q)\dif\mu_{\BS x,t}^{\prime q,\omega}(\BS z)
\dif\mathbb P(\omega)\dif\rho(q),
\end{equation}
\end{subequations}
where $\mu'$ is the derivative of $\mu$ with respect to the parameter
$\BS x$. Observe that now the solution of
\eqref{eq:dyn_sys_nonauto_slow_limiting_x}, starting from a point $\BS
z$, is given by the cocycle $\phi$:
\begin{equation}
\BS z(s)=\phi_{\BS x,t}^{q,\omega,s}\BS z,
\end{equation}
where $q$ and $\omega$ denote parametric dependence on the starting
time and probability space outcome for the Wiener process $\BS
W_s(\omega)$, which means that $\phi$ satisfies
\begin{equation}
\dif\phi_{\BS x,t}^{q,\omega,s}\BS z=\BS{\tilde G}(\BS x,\phi_{\BS x,t}
^{q,\omega,s}\BS z,t,\theta^sq)\dif s+\BS{\tilde\sigma}(\phi_{\BS x,t}
^{q,\omega,s}\BS z,t,\theta^sq)\dif\BS W_s(\omega). 
\end{equation}
Now, let us denote the tangent map of the cocycle $\phi$ as
\begin{equation}
\BS{\tilde T}_{\BS x,t,\BS z}^{q,\omega,s}=\parderiv{}{\BS z}
\phi_{\BS x,t}^{q,\omega,s}\BS z,
\end{equation}
which obviously obeys the chain rule
\begin{equation}
\BS{\tilde T}_{\BS x,t,\BS z}^{q,\omega,s+r}=\BS{\tilde T}_{\BS
  x,t,\phi_{\BS x,t}^{q,\omega,r}\BS z}
^{\theta^r q,\vartheta^r \omega,s}\;\;\BS{\tilde T}_{\BS
  x,t,\BS z}^{q,\omega,r},\quad\forall(q,\omega)\in(\mathcal
Q,\Omega)
\end{equation}
and satisfies the equation
\begin{equation}
\dif\BS{\tilde T}_{\BS x,t,\BS z}^{q,\omega,s}=\left[\parderiv{\BS{\tilde G}}
{\BS z}(\BS x,\phi_{\BS x,t}^{q,\omega,s}\BS z,t,\theta^sq)\dif s
+\parderiv{\BS{\tilde\sigma}}{\BS z}(\phi_{\BS x,t}
^{q,\omega,s}\BS z,t,\theta^sq)\dif\BS W_s(\omega)\right]
\BS{\tilde T}_{\BS x,t,\BS z}^{q,\omega,s}. 
\end{equation}
Then, through the Duhamel's principle, the $\BS x$-derivative of $\phi$
after elapsed time $a$ is computed as
\begin{equation}
\label{eq:phi_x}
\parderiv{}{\BS x}\phi_{\BS x,t}^{q,\omega,a}\BS z=\int_0^a
\BS{\tilde T}_{\BS x,t,\phi_{\BS x,t}^{q,\omega,s}\BS z}^{\theta^s q,
\vartheta^s\omega,a-s}\parderiv{\BS{\tilde G}}{\BS x}(\BS x,
\phi_{\BS x,t}^{q,\omega,s}\BS z,t,\theta^s q)\dif s,
\end{equation}
while the average linear response of $\BS F$ to the changes in
$\mu_{\BS x,t}^{q,\omega}$, which is after the elapsed time $a$ since
$\BS x$ was changed, is given, in the pullback sense, as
\begin{equation}
\label{eq:ext_temp1}
\begin{split}
\int_{\mathcal A_{\BS x,t}(q,\omega)}\BS{\tilde F}(\BS x,\BS z,t,q)\dif
\mu_{\BS x,t}^{\prime q,\omega}(\BS z)=\int_0^a\int_{\mathcal A_{\BS x,t}
(\theta^{-a}q,\vartheta^{-a}\omega)}\parderiv{\BS{\tilde F}}{\BS z}
(\BS x,\BS z,t,q)\BS{\tilde T}_{\BS x,t,\phi_{\BS x,t} ^{q,\omega,s-a}\BS z}
^{\theta^{s-a} q,\vartheta^{s-a}\omega,a-s} \times\\\times
\parderiv{\BS{\tilde G}}{\BS x}(\BS x,\phi_{\BS x,t} ^{q,\omega,s-a}
\BS z,t,\theta^{s-a}q) \dif\mu_{\BS x,t}^{\theta^{-a}q,\vartheta^{-a}\omega}
(\BS z)\dif s.
\end{split}
\end{equation}
After applying \eqref{eq:mu_inv}, the above formula becomes
\begin{equation}
\begin{split}
\int_{\mathcal A_{\BS x,t}(q,\omega)}\BS{\tilde F}(\BS x,\BS z,t,q)
\dif\mu_{\BS x,t}^{\prime q,\omega}(\BS z)=\int_0^a\int_{\mathcal
A_{\BS x,t}(\theta^{s-a}q,\vartheta^{s-a}\omega)}\parderiv{\BS{\tilde F}}
{\BS z}(\BS x,\phi_{\BS x,t}^{\theta^{s-a}q,\theta^{s-a} \omega,a-s}\BS
z,t,q)\times\\\times\BS{\tilde T}_{\BS x,t,\BS z}^{\theta^{s-a}q,
\vartheta^{s-a} \omega,a-s}\parderiv{\BS{\tilde G}}{\BS x}(\BS x,
\BS z,t,\theta^{s-a} q)\dif\mu_{\BS x,t}^{\theta^{s-a}q,\vartheta^{s-a}
\omega} (\BS z)\dif s.
\end{split}
\end{equation}
The further average over $\rho$ and $\mathbb P$ yields
\begin{equation}
\begin{split}
\int_{\mathcal Q}\int_\Omega&\int_{\mathcal A_{\BS x,t}(q,\omega)}
\BS{\tilde F}(\BS x,\BS z,t,q)\dif\mu_{\BS x,t}^{\prime q,\omega}
(\BS z)\dif\mathbb P(\omega)\dif\rho(q)=\\=&\int_0^a
\int_{\mathcal Q}\int_\Omega\int_{\mathcal
  A_{\BS x,t}(\theta^{s-a}q,\vartheta^{s-a}\omega)}
\parderiv{\BS{\tilde F}}{\BS z}(\BS x,\phi_{\BS x,t}
^{\theta^{s-a}q,\vartheta^{s-a}\omega,a-s}\BS z,t,q) \BS{\tilde T}_{\BS
  x,t,\BS z}^{\theta^{s-a}q,\vartheta^{s-a}\omega,a-s}
\times\\&\times\parderiv{\BS{\tilde G}}{\BS x}(\BS x,\BS
z,t,\theta^{s-a}q) \dif\mu_{\BS
  x,t}^{\theta^{s-a}q,\vartheta^{s-a}\omega}(\BS z)
\dif\mathbb P(\omega)\dif\rho(q)\dif s,
\end{split}
\end{equation}
where, after using the invariance of $\rho$ and $\mathbb P$ under
$\theta$ and $\vartheta$, respectively, and replacing $a-s$ with $s$
in the integral over $s$, finally obtain
\begin{equation}
\begin{split}
\int_{\mathcal Q}\int_\Omega\int_{\mathcal A_{\BS x,t}(q,\omega)}
\BS{\tilde F}(\BS x,\BS z,t,q)\dif\mu_{\BS x,t}^{\prime q,\omega}
(\BS z)\dif\mathbb P(\omega)\dif\rho(q)=\\=\int_0^a\int_{\mathcal Q}
\int_\Omega\int_{\mathcal A_{\BS x,t}(q,\omega)}\parderiv{}{\BS z}
\BS{\tilde F}(\BS x,\phi_{\BS x,t}^{q,\omega,s}\BS z,t,\theta^sq)
\times\\\times\BS{\tilde T}_{\BS x,t,\BS z}^{q,\omega,s}\parderiv{}
{\BS x}\BS{\tilde G}(\BS x,\BS z,t,q)\dif\mu_{\BS x,t}^{q,\omega}
(\BS z)\dif\mathbb P(\omega)\dif\rho(q)\dif s.
\end{split}
\end{equation}
Then the term $\BS H(\BS x,t)$ in
\eqref{eq:dyn_sys_nonauto_slow_tangent} is formally given by sending
the response time $a$ in the formula above to infinity:
\begin{equation}
\begin{split}
\BS H(\BS x,t)=\int_0^\infty\int_{\mathcal Q}\int_\Omega
\int_{\mathcal A_{\BS x,t}(q,\omega)}\parderiv{}{\BS z} \BS{\tilde
  F}(\BS x,\phi_{\BS x,t}^{q,\omega,s}\BS z,t,\theta^sq)
\times\\\times\BS{\tilde T}_{\BS x,t,\BS z}^{q,\omega,s}\parderiv{}
{\BS x}\BS{\tilde G}(\BS x,\BS z,t,q)\dif\mu_{\BS x,t}^{q,\omega}(\BS z)
\dif\mathbb P(\omega)\dif\rho(q)\dif s.
\end{split}
\end{equation}
The assumption of the joint ergodicity of $\mathcal Q$ and $\Omega$ is
that a sufficiently long trajectory $s$ completely samples $\mathcal
Q$ and $\Omega$ through $\theta^s$ and $\vartheta^s$, respectively,
with the corresponding statistical weights given by the product of the
measures $\rho$ and $\mathbb P$. Under this assumption, one can
replace the measure averages with the time averages, leading to
\eqref{eq:explicit_time_averaged}.

\end{document}